\documentclass{svjour3}

\usepackage{mathptmx}

\usepackage{amssymb}
\usepackage{amsmath}
\usepackage[colorlinks, linkcolor=blue, citecolor=blue, pdfborder={0 0 0 [0 0]}]{hyperref}

\makeatletter
\newcommand{\mylabel}[2]{#2\def\@currentlabel{#2}\label{#1}}
\makeatother

\newcommand{\Lip}{\operatorname{Lip}}
\newcommand{\rd}{\mathrm{d}}
\newcommand{\dist}{\operatorname{dist}}

\begin{document}

\title{Path-Dependent Hamilton--Jacobi Equations: The Minimax Solutions Revised
\thanks{The work was performed as part of research conducted in the Ural Mathematical Center.}}

\author{Mikhail I. Gomoyunov \and Nikolai Yu. Lukoyanov \and Anton R. Plaksin}

\authorrunning{M. I. Gomoyunov, N. Yu. Lukoyanov, A. R. Plaksin}

\institute{M. I. Gomoyunov, N. Yu. Lukoyanov, A. R. Plaksin \at
    Krasovskii Institute of Mathematics and Mechanics of the Ural Branch of the Russian Academy of Sciences,
    S. Kovalevskaya Str., 16, Ekaterinburg, Russia \\
    Ural Federal University,
    Mira Str., 32, Ekaterinburg, Russia \\
    \email{m.i.gomoyunov@gmail.com, nyul@imm.uran.ru, a.r.plaksin@gmail.com}}

\date{Received: date / Accepted: date}

\maketitle

\begin{abstract}
    Motivated by optimal control problems and differential games for functional differential equations of retarded type, the paper deals with a Cauchy problem for a path-dependent Hamilton--Jacobi equation with a right-end boundary condition.
    Minimax solutions of this problem are studied.
    The existence and uniqueness result is obtained under assumptions that are weaker than those considered earlier.
    In contrast to previous works, on the one hand, we do not require any properties concerning positive homogeneity of the Hamiltonian in the impulse variable, and on the other hand, we suppose that the Hamiltonian satisfies a Lipshitz continuity condition with respect to the path variable in the uniform (supremum) norm.
    The progress is related to the fact that a suitable Lyapunov--Krasovskii functional is built that allows to prove a comparison principle.
    This functional is in some sense equivalent to the square of the uniform norm of the path variable and, at the same time, it possesses appropriate smoothness properties.
    In addition, the paper provides non-local and infinitesimal criteria of minimax solutions, their stability with respect to perturbations of the Hamiltonian and the boundary functional, as well as consistency of the approach with the non-path-dependent case.
    Connection of the problem statement under consideration with some other possible statements (regarding the choice of path spaces and derivatives used) known in the theory of path-dependent Hamilton--Jacobi equations is discussed.
    Some remarks concerning viscosity solutions of the studied Cauchy problem are given.


\end{abstract}

\section{Introduction}
    In the paper, for a functional possessing a non-anticipativity property, we consider a path-dependent Hamilton--Jacobi equation with co-invariant derivatives and study a Cauchy problem for this equation under a right-end boundary condition.

    Co-invariant derivatives of non-anticipative functionals were initially introduced in \cite{Kim_1985} for the purposes of the stability theory of functional differential equations of retarded type (see, e.g., \cite{Krasovskii_1963}).
    Later, these derivatives found application in other parts of the theory of functional differential equations (see, e.g., \cite{Kim_1999} and the references therein).
    Among similar notions of derivatives of non-anticipative functionals, let us note Clio derivatives \cite{Aubin_Haddad_2002} and horizontal and vertical derivatives \cite{Dupire_2009}.
    On the other hand, it should be noted also that the derivatives of this kind are inherently different from classical Fr\'{e}chet derivatives of functionals (in this regard, see, e.g., \cite{Ji_Yang_2015}).

    Examples of Hamilton--Jacobi equations considered in the paper are the Bellman equations and the Bellman--Isaacs equations associated respectively to optimal control problems and differential games for functional differential equations of retarded type (see, e.g., \cite{Kim_1999,Lukoyanov_2000_JAMM,Lukoyanov_2003_2,Kaise_2015,Plaksin_2019_IFAC} and also \cite{Aubin_Haddad_2002,Bayraktar_Keller_2018,Zhou_2019}).
    Since the value and the Hamiltonian in these problems are non-anticipative functionals depending on a history of a path, the corresponding Hamil\-ton--Jacobi equations are often called non-anticipative or path-dependent.
    In cases when the value functional has some additional smoothness properties (namely, it is co-invariantly smooth), it is characterized as a classical solution of the Cauchy problem for the associated path-dependent Hamilton--Jacobi equation and the natural right-end boundary condition.
    However, as a rule, the value functional does not possess the required smoothness properties, it can not be treated as a classical solution of this Cauchy problem, and, thus, appropriate generalized solutions should be considered.
    Let us mention also that path-dependent Hamilton--Jacobi equations arise in optimal control problems and differential games for ordinary differential equations but with path-dependent cost (see, e.g., \cite{Kaise_Kato_Takahashi_2018} and also \cite{Krasovskii_Krasovskii_1995,Lukoyanov_Gomoyunov_2019_DGAA} for related issues).

    In the theory of Hamilton--Jacobi equations, various notions of generalized solutions are proposed, including minimax (see, e.g., \cite{Subbotin_1991_NA,Subbotin_1995}) and viscosity (see, e.g., \cite{Crandall_Lions_1983,Crandall_Evans_Lions_1984}) solutions.
    The present paper deals with the minimax solutions.
    The theory of minimax solutions of Hamilton--Jacobi equations originates in unification constructions of the positional differential games (see, e.g., \cite{Krasovskii_1976,Krasovskii_Subbotin_1988}).
    The term ``minimax'' refers to the fact that minimum and maximum operations play an important role in the basic elements of this theory.

    The theory of minimax solutions of Hamilton--Jacobi equations with partial derivatives is presented in \cite{Subbotin_1991_Eng,Subbotin_1995} (see also \cite{Subbotin_1991_NA}).
    Minimax solutions of Cauchy problems for path-dependent Hamilton--Jacobi equations with co-invariant derivatives are studied in, e.g., \cite{Krasovskii_Lukoyanov_2000_IMM_Eng,Lukoyanov_2001_DE_Eng,Lukoyanov_2003_1,Lukoyanov_2010_IMM_Eng_1} (see also \cite{Lukoyanov_2011_Eng}) under various groups of assumptions.
    Firstly, following \cite{Subbotin_1991_Eng} (see also \cite{Subbotin_1991_NA}), the cases were considered when the Hamiltonian is homogeneous in the impulse variable \cite{Lukoyanov_2001_DE_Eng} or when it is not homogenous but the problem can be reduced to an auxiliary problem with a homogeneous Hamiltonian \cite{Krasovskii_Lukoyanov_2000_IMM_Eng}.
    After that, following \cite{Subbotin_1995}, a fully non-homogeneous case was investigated \cite{Lukoyanov_2003_1,Lukoyanov_2010_IMM_Eng_1}.
    However, assumed in \cite{Lukoyanov_2003_1,Lukoyanov_2010_IMM_Eng_1} Lipschitz continuity properties of the Hamiltonian with respect to the path variable, in the context of applications to optimal control problems and differential games, allow to consider distributed and constant concentrated delays only.
    These studies are mainly focused on well-posedeness of minimax solutions (existence, uniqueness, and stability with respect to perturbations of the Hamiltonian and the boundary functional), consistency of minimax solutions with classical solutions, as well as non-local and infinitesimal criteria of minimax solutions.
    In, e.g., \cite{Lukoyanov_2000_JAMM,Lukoyanov_2001_PMM_Eng,Lukoyanov_2003_2,Lukoyanov_2004_PMM_Eng,Lukoyanov_2006_IMM_Eng,Lukoyanov_2010_IMM_Eng_1}, these results are applied to characterize the game value functional and construct optimal positional (closed-loop) strategies in zero-sum differential games for functional differential equations of retarded type.
    Concerning viscosity solutions of Cauchy problems for path-dependent Hamilton--Jacobi equations with co-invariant derivatives, we refer to \cite{Lukoyanov_2007_IMM_Eng,Lukoyanov_2010_IMM_Eng_1}, where also a connection between viscosity and minimax solutions of such problems is investigated.

    Let us note that approximations of minimax solutions of path-dependent Hamil\-ton--Ja\-cobi equations with co-invariant derivatives by means of minimax solutions of auxiliary Hamilton--Jacobi equations with partial derivatives are considered in \cite{Krasovskii_Lukoyanov_2000_IMM_Eng,Lukoyanov_2011_Eng,Gomoyunov_Lukoyanov_Plaksin_2019_IMM_Eng}.
    In \cite{Bayraktar_Keller_2018}, minimax solutions of path-dependent Hamilton--Jacobi equations in an infinite dimensional setting are studied.
    In \cite{Lukoyanov_Gomoyunov_Plaksin_2017_Doklady,Plaksin_2019_DE_Eng}, the theory of minimax solutions is developed for path-dependent Hamilton--Jacobi equations arising in optimal control problems and differential games for functional differential equations of neutral type.

    The present paper generalizes the results \cite{Krasovskii_Lukoyanov_2000_IMM_Eng,Lukoyanov_2001_DE_Eng,Lukoyanov_2003_1,Lukoyanov_2010_IMM_Eng_1} to the case of weaker assumptions on the Hamiltonian.
    Namely, in contrast to \cite{Krasovskii_Lukoyanov_2000_IMM_Eng,Lukoyanov_2001_DE_Eng}, any additional suppositions concerning homogeneity of the Hamiltonian are not required.
    Compared to \cite{Lukoyanov_2003_1,Lukoyanov_2010_IMM_Eng_1}, it is assumed that the Hamiltonian satisfies a Lipshitz continuity condition with respect to the path variable in the uniform (supremum) norm, which, in particular, allows to cover not only the cases of distributed and constant concentrated delays, but also the case of time-varying delays.
    The main result of the paper is an existence and uniqueness theorem for minimax solutions of the considered Cauchy problem.
    The progress is related to the fact that a suitable Lyapunov--Krasovskii functional (see, e.g., \cite{Lukoyanov_2004_PMM_Eng,Lukoyanov_2006_IMM_Eng}) is found, which is in some sense equivalent to the square of the uniform norm and, at the same time, possesses appropriate smoothness properties.
    This functional is built on the basis of the constructions from \cite[Sect.~7.5]{Subbotin_1995} and \cite{Zhou_2020_1}.

    The paper is organized as follows.
    In Sect.~\ref{section_problem_statement}, we recall definitions of non-anticipative functionals and their co-invariant derivatives, formulate a Cauchy problem for a path-de\-pend\-ent Hamilton--Jacobi equation with co-invariant derivatives, which is the subject of the paper, and describe the basic assumptions.
    As a motivation, we give an optimal control problem for functional differential equations of retarded type that leads to the considered Hamilton--Jacobi equations.
    In Sect.~\ref{section_minimax_solution}, we define minimax solutions of the studied Cauchy problem, provide non-local and infinitesimal criteria of such solutions, and discuss consistency of minimax solutions with classical solutions.
    Sect.~\ref{section_well_posedeness} is devoted to well-posedeness of minimax solutions.
    We prove the corresponding existence and uniqueness theorem, which, in particular, allows to establish stability of minimax solutions with respect to perturbations of the Hamiltonian and the boundary functional.
    In the end of this section, we show that, in the non-path-dependent case, the considered (path-dependent) minimax solutions are consistent with minimax solutions of the corresponding Hamilton--Jacobi equations with partial derivatives.
    In Sect.~\ref{section_remarks_on_problem_statement}, we discuss connection of the problem statement considered in the paper with some other possible approaches (regarding the choice of path spaces and derivatives used) developed in the theory of path-dependent Hamilton--Jacobi equations.
    In Sect.~\ref{section_remarks_homogenous}, we compare the results of the paper with those obtained in \cite{Krasovskii_Lukoyanov_2000_IMM_Eng,Lukoyanov_2001_DE_Eng}.
    In Sect.~\ref{section_viscosity}, we give some remarks concerning viscosity solutions of the studied Cauchy problem.

    Summarizing, we conclude that the results of the paper allow to obtain the theory of minimax solutions of path-dependent Hamilton--Jacobi equations with co-invariant derivatives that is most consistent with the theory of minimax solutions of Hamilton--Jacobi equations with partial derivatives.

\section{Problem Statement}
\label{section_problem_statement}

    Let us introduce some notations and terminology used in the paper.
    Let $n \in \mathbb{N}$ be fixed, and let $\mathbb{R}^n$ be the Euclidean space of $n$-dimensional vectors with the inner product $\langle \cdot, \cdot \rangle$ and the norm $\|\cdot\|$.
    Given $t_1$, $t_2 \in \mathbb{R}$, $t_1 \leq t_2$, for a function $x: [t_1, t_2] \to \mathbb{R}^n$, we use the notation $x(\cdot)$ if this function is considered as an element of some functional space, while the value of this function at a particular point $t \in [t_1, t_2]$ is denoted, as usual, by $x(t)$.
    Further, let $C([t_1, t_2], \mathbb{R}^n)$ be the Banach space of continuous functions $x: [t_1, t_2] \to \mathbb{R}^n$ with the uniform norm
    \begin{equation} \label{norm_C}
        \|x(\cdot)\|_{[t_1, t_2]}
        = \max_{t \in [t_1, t_2]} \|x(t)\|,
        \quad x(\cdot) \in C([t_1, t_2], \mathbb{R}^n).
    \end{equation}
    Finally, let $\dot{x}(t)$ stand for the derivative $\rd x(t) / \rd t$ of a function $x: [t_1, t_2] \to \mathbb{R}^n$ at a point $t \in (t_1, t_2)$.

    \subsection{Non-Anticipative Functionals}

        Let $h > 0$ and $T > 0$ be fixed.
        Consider the set $[0, T] \times C([- h, T], \mathbb{R}^n)$ with the standard product metric
        \begin{equation} \label{dist}
            \dist\big( (t, x(\cdot)), (\tau, y(\cdot)) \big)
            = |t - \tau| + \|x(\cdot) - y(\cdot)\|_{[- h, T]}.
        \end{equation}

        A functional $\varphi: [0, T] \times C([- h, T], \mathbb{R}^n) \to \mathbb{R}$ is said to be {\it non-anticipative} if, for every $t \in [0, T)$ and $x(\cdot)$, $y(\cdot) \in C([- h, T], \mathbb{R}^n)$, the equality
        \begin{equation} \label{x_equals_y}
            x(\tau) = y(\tau),
            \quad \tau \in [- h, t],
        \end{equation}
        implies that $\varphi(t, x(\cdot)) = \varphi(t, y(\cdot))$.

        Some remarks about continuity properties of non-anticipative functionals are given in Sect.~\ref{subsection_continuity} below.

    \subsection{Co-Invariant Derivatives}

        For every $(t, x(\cdot)) \in [0, T] \times C([- h, T], \mathbb{R}^n)$, consider the set $\Lip(t, x(\cdot))$ consisting of functions $y(\cdot) \in C([- h, T], \mathbb{R}^n)$ that satisfy \eqref{x_equals_y} and are Lipschitz continuous on $[t, T]$.

        In accordance with, e.g., \cite[Definition~2.4.1]{Kim_1999} and \cite{Krasovskii_Lukoyanov_2000_IMM_Eng,Lukoyanov_2000_JAMM}, we say that a functional $\varphi: [0, T] \times C([- h, T], \mathbb{R}^n) \to \mathbb{R}$ is {\it co-invariantly} ($ci$-) {\it differentiable} at a point $(t, x(\cdot)) \in [0, T) \times C([- h, T], \mathbb{R}^n)$ if there exist $\partial_t \varphi(t, x(\cdot)) \in \mathbb{R}$ and $\nabla \varphi(t, x(\cdot)) \in \mathbb{R}^n$ such that, for every $y(\cdot) \in \Lip(t, x(\cdot))$, the following relation is valid for any $\tau \in (t, T]$:
        \begin{equation} \label{ci_derivatives_definition}
            \varphi(\tau, y(\cdot)) - \varphi(t, x(\cdot))
            = \partial_t \varphi(t, x(\cdot)) (\tau - t) + \langle \nabla \varphi(t, x(\cdot)), y(\tau) - x(t) \rangle +  o(\tau - t),
        \end{equation}
        where the function $o(\delta) \in \mathbb{R}$, $\delta > 0$, may depend on $y(\cdot)$ and $o(\delta) / \delta \to 0$ as $\delta \downarrow 0$.
        In this case, the values $\partial_t \varphi(t, x(\cdot))$ and $\nabla \varphi(t, x(\cdot))$ are called $ci$-derivatives of $\varphi$ at $(t, x(\cdot))$.

        Note that, if  a functional $\varphi: [0, T] \times C([- h, T], \mathbb{R}^n) \to \mathbb{R}$ is $ci$-differentiable at some points $(t, x(\cdot))$, $(t, y(\cdot)) \in [0, T) \times C([- h, T], \mathbb{R}^n)$ satisfying \eqref{x_equals_y}, then the equalities
        \begin{equation*}
            \varphi(t, x(\cdot))
            = \varphi(t, y(\cdot)),
            \quad \partial_t \varphi(t, x(\cdot))
            = \partial_t \varphi(t, y(\cdot)),
            \quad \nabla \varphi(t, x(\cdot))
            = \nabla \varphi(t, y(\cdot))
        \end{equation*}
        hold.
        Therefore, if $\varphi$ is $ci$-differentiable at every point $(t, x(\cdot)) \in [0, T) \times C([- h, T], \mathbb{R}^n)$, then the functional $\varphi$ itself and its $ci$-derivatives
        \begin{equation} \label{ci_derivatives}
            \begin{array}{c}
                [0, T) \times C([- h, T], \mathbb{R}^n) \ni (t, x(\cdot)) \mapsto \partial_t \varphi (t, x(\cdot)) \in \mathbb{R}, \\[0.5em]
                [0, T) \times C([- h, T], \mathbb{R}^n) \ni (t, x(\cdot)) \mapsto \nabla \varphi (t, x(\cdot)) \in \mathbb{R}^n
            \end{array}
        \end{equation}
        are automatically non-anticipative.

        We say that a functional $\varphi: [0, T] \times C([- h, T], \mathbb{R}^n) \to \mathbb{R}$ is {\it $ci$-smooth} if it is continuous, $ci$-differentiable at every point $(t, x(\cdot)) \in [0, T) \times C([- h, T], \mathbb{R}^n)$, and its $ci$-derivatives \eqref{ci_derivatives} are continuous.
        The following proposition gives one of the key properties of $ci$-smooth functionals (see, e.g., \cite[Theorem~7.1.1]{Kim_1999}, \cite[Lemma~2.1]{Lukoyanov_2011_Eng}, and also the proofs of \cite[Lemma~3]{Plaksin_2019_IFAC} and \cite[Lemma~9.2]{Gomoyunov_2019_SIAM}).
        \begin{proposition} \label{proposition_ci_smooth}
            If a functional $\varphi: [0, T] \times C([- h, T], \mathbb{R}^n) \to \mathbb{R}$ is $ci$-smooth, then, for every $(t, x(\cdot)) \in [0, T) \times C([- h, T], \mathbb{R}^n)$ and $y(\cdot) \in \Lip(t, x(\cdot))$, the equality below holds:
            \begin{equation} \label{proposition_ci_smooth_main}
                \varphi(\tau, y(\cdot))
                = \varphi(t, x(\cdot))
                + \int_{t}^{\tau} \partial_t \varphi(\xi, y(\cdot)) \, \rd \xi
                + \int_{t}^{\tau} \langle \nabla \varphi(\xi, y(\cdot)), \dot{y}(\xi) \rangle \, \rd \xi,
                \quad \tau \in [t, T).
            \end{equation}
        \end{proposition}

        Some additional remarks concerning the notion of $ci$-differentiability and its connection with some other approaches to differentiation of non-anticipative functionals are given in Sect.~\ref{subsection_derivatives} below.

    \subsection{Path-Dependent Hamilton--Jacobi Equation}

        Let a {\it Hamiltonian} $H: [0, T] \times C([- h, T], \mathbb{R}^n) \times \mathbb{R}^n \to \mathbb{R}$ and a {\it boundary functional} $\sigma: C([- h, T], \mathbb{R}^n) \to \mathbb{R}$ be given, and suppose that, for every fixed $s \in \mathbb{R}^n$, the functional
        \begin{equation} \label{H_s}
            [0, T] \times C([-h, T], \mathbb{R}^n) \ni (t, x(\cdot)) \mapsto H(t, x(\cdot), s) \in \mathbb{R}
        \end{equation}
        is non-anticipative.
        In the paper, we study the {\it Cauchy problem} for the {\it Hamilton--Jacobi equation} with $ci$-derivatives
        \begin{equation} \label{HJ}
            \partial_t \varphi(t, x(\cdot)) + H \big( t, x(\cdot), \nabla \varphi(t, x(\cdot)) \big)
            = 0,
            \quad (t, x(\cdot)) \in [0, T) \times C([- h, T], \mathbb{R}^n),
        \end{equation}
        and the right-end {\it boundary condition}
        \begin{equation} \label{boundary_condition}
            \varphi(T, x(\cdot))
            = \sigma(x(\cdot)),
            \quad x(\cdot) \in C([- h, T], \mathbb{R}^n).
        \end{equation}
        The unknown here is a non-anticipative functional $\varphi: [0, T] \times C([- h, T], \mathbb{R}^n) \to \mathbb{R}$.

        An example of \eqref{HJ} is the Bellman equation associated to the optimal control problem described by an initial data $(t, x(\cdot)) \in [0, T] \times C([- h, T], \mathbb{R}^n)$, the dynamic equation
        \begin{equation} \label{system}
            \dot{y}(\tau)
            = f(\tau, y(\cdot), u(\tau)),
            \quad \tau \in [t, T],
        \end{equation}
        with initial condition \eqref{x_equals_y} specified by $(t, x(\cdot))$, and the cost functional
        \begin{equation} \label{cost_functional}
            J
            = \sigma(y(\cdot)) - \int_{t}^{T} g(\tau, y(\cdot), u(\tau)) \, \rd \tau
        \end{equation}
        to be minimized.
        Here, $\tau$ is time; $y(\tau) \in \mathbb{R}^n$ is the current state; $u(\tau) \in U$ is the current control, $U \subset \mathbb{R}^m$ is a compact set, $m \in \mathbb{N}$.
        The set of admissible controls $\mathcal{U}$ consists of all (Lebesgue) measurable functions $u: [0, T] \to U$.
        By a path of the considered dynamical system, we mean a function $y(\cdot) \in \Lip(t, x(\cdot))$ that together with a control $u(\cdot) \in \mathcal{U}$ satisfies equation \eqref{system} for almost every (a.e.) $\tau \in [t, T]$.

        It is assumed that the following conditions are fulfilled:
        \begin{description}
            \item[\rm{(\mylabel{A.1}{$A.1$})}]
                The functions $f: [0, T] \times C([-h, T], \mathbb{R}^n) \times U \to \mathbb{R}^n$, $g: [0, T] \times C([-h, T], \mathbb{R}^n) \times U \to \mathbb{R}$, and $\sigma: C([-h, T], \mathbb{R}^n) \to \mathbb{R}$ are continuous.

            \item[\rm{(\mylabel{A.2}{$A.2$})}]
                There exists a constant $c > 0$ such that
                \begin{equation*}
                    \|f(\tau, y(\cdot), u)\|
                    \leq c (1 + \max_{\xi \in [- h, \tau]} \|y(\xi)\|)
                \end{equation*}
                for every $\tau \in [0, T]$, $y(\cdot) \in C([- h, T], \mathbb{R}^n)$, and $u \in U$.

            \item[\rm{(\mylabel{A.3}{$A.3$})}]
                For any compact set $D \subset C([- h, T], \mathbb{R}^n)$, there exists a number $\lambda > 0$ such that
                \begin{equation*}
                    \|f(\tau, y(\cdot), u) - f(\tau, z(\cdot), u)\| + |g(\tau, y(\cdot), u) - g(\tau, z(\cdot), u)|
                    \leq \lambda \max_{\xi \in [- h, \tau]} \|y(\xi) - z(\xi)\|
                \end{equation*}
                for every $\tau \in [0, T]$, $y(\cdot)$, $z(\cdot) \in D$, and $u \in U$.
        \end{description}

        Note that condition \eqref{A.3} implies non-anticipativity of the functions
        \begin{equation*}
            \begin{array}{c}
                [0, T] \times C([- h, T], \mathbb{R}^n) \ni (\tau, y(\cdot)) \mapsto f(\tau, y(\cdot), u) \in \mathbb{R}^n, \\[0.5em]
                [0, T] \times C([- h, T], \mathbb{R}^n) \ni (\tau, y(\cdot)) \mapsto g(\tau, y(\cdot), u) \in \mathbb{R}
            \end{array}
        \end{equation*}
        for every fixed $u \in U$.
        Thus, dynamic equation \eqref{system} is actually a functional differential equation of retarded type.
        Observe that conditions \eqref{A.1} to \eqref{A.3} are quite natural and allow to cover such important for applications cases as {\it constant concentrated delays}, when, e.g., the function $f$ (similarly, the function $g$) can be represented as
        \begin{equation} \label{f_1}
            f(\tau, y(\cdot), u)
            = f_1 \big(\tau, y(\tau), y(\tau - h), u \big),
        \end{equation}
        {\it time-varying concentrated delays}, when, e.g.,
        \begin{equation} \label{f_2}
            f(\tau, y(\cdot), u)
            = f_2 \big(\tau, y(\tau), y(\tau - k(\tau)), u \big),
            \quad 0 < k(\tau) \leq h,
        \end{equation}
        and {\it distributed delays}, when, e.g.,
        \begin{equation}  \label{f_3}
            f(\tau, y(\cdot), u)
            = f_3 \Big(\tau, y(\tau), \int_{- h}^\tau K(\tau, \xi, y(\xi)) \, \rd \xi, u \Big).
        \end{equation}

        For every $(t, x(\cdot)) \in [0, T] \times C([- h, T], \mathbb{R}^n)$ and $u(\cdot) \in \mathcal{U}$, conditions \eqref{A.1} to \eqref{A.3} provide existence and uniqueness of a path $y(\cdot)$ of the system.
        Let $J(t, x(\cdot), u(\cdot))$ denote the corresponding value of cost functional \eqref{cost_functional}.
        Then, the optimal result functional (the value functional) is given by
        \begin{equation*}
            \varphi(t, x(\cdot))
            = \inf_{u(\cdot) \in \mathcal{U}} J(t, x(\cdot), u(\cdot)),
            \quad (t, x(\cdot)) \in [0, T] \times C([- h, T], \mathbb{R}^n).
        \end{equation*}
        By construction, this functional $\varphi$ is non-anticipative.
        For $t = T$, it satisfies boundary condition \eqref{boundary_condition}.
        For $t \in [0, T)$, it is characterized (see, e.g., \cite[Sect.~10.2]{Kim_1999} and \cite[Theorem~3.2]{Zhou_2019}) by the following relation, expressing the dynamic programming principle in the considered optimal control problem:
        \begin{equation} \label{DPP}
            \varphi(t, x(\cdot))
            = \inf_{u(\cdot) \in \mathcal{U}} \Big( \varphi(\tau, y(\cdot)) - \int_{t}^{\tau} g(\xi, y(\cdot), u(\xi)) \, \rd \xi \Big),
            \quad \tau \in (t, T].
        \end{equation}

        Now, if we assume that the value functional $\varphi$ is $ci$-smooth, then, based on \eqref{proposition_ci_smooth_main} and \eqref{DPP}, we derive (see, e.g., \cite[Sect.~10.2]{Kim_1999} and \cite[Theorem~3.3]{Zhou_2020_1}) that it satisfies equation \eqref{HJ} with the Hamiltonian
        \begin{equation} \label{Hamiltonian_OCP}
            H(t, x(\cdot), s)
            = \min_{u \in U} \big( \langle s, f(t, x(\cdot), u) \rangle - g(t, x(\cdot), u) \big).
        \end{equation}
        Equation \eqref{HJ} with Hamiltonian \eqref{Hamiltonian_OCP} is called the Hamilton--Jacobi--Bellman equation or just the Bellman equation associated to optimal control problem \eqref{system}, \eqref{cost_functional}.

        Let us observe that, in differential games theory for functional differential equations of retarded type, equation \eqref{HJ} arises as the Bellman--Isaacs equation (see, e.g., \cite{Lukoyanov_2001_PMM_Eng,Lukoyanov_2003_2,Lukoyanov_2010_IMM_Eng_1} and also \cite{Lukoyanov_2011_Eng}).

        Let us emphasize that, in problem \eqref{HJ}, \eqref{boundary_condition}, the Hamiltonian $H$, the $ci$-derivatives $\partial_t \varphi$ and $\nabla \varphi$, as well as the sought functional $\varphi$ itself depend at any fixed $t \in [0, T)$ on a restriction $x_t(\cdot)$ of $x(\cdot)$ to $[- h, t]$ given by
        \begin{equation} \label{x_t}
            x_t(\tau)
            = x(\tau),
            \quad \tau \in [- h, t].
        \end{equation}
        In relation to optimal control problems and differential games, it means that these values depend on a ``path segment'' $y_t(\cdot) = x_t(\cdot)$, which can be understood as history until time $t$ of a path $y(\cdot)$.
        In this context, equation \eqref{HJ} is called {\it non-anticipative} or {\it path-dependent} Hamilton--Jacobi equation.

    \subsection{Assumptions}

        Throughout the paper, we assume that the Hamiltonian $H$ and the boundary functional $\sigma$ in problem \eqref{HJ}, \eqref{boundary_condition} satisfy the following conditions:
        \begin{description}
            \item[\rm{(\mylabel{B.1}{$B.1$})}]
                The functionals $H$ and $\sigma$ are continuous.

            \item[\rm{(\mylabel{B.2}{$B.2$})}]
                There exists a constant $c > 0$ such that
                \begin{equation*}
                    |H(t, x(\cdot), s) - H(t, x(\cdot), r)|
                    \leq c (1 + \max_{\tau \in [- h, t]} \|x(\tau)\|) \|s - r\|
                \end{equation*}
                for every $t \in [0, T]$, $x(\cdot) \in C([- h, T], \mathbb{R}^n)$, and $s$, $r \in \mathbb{R}^n$.

            \item[\rm{(\mylabel{B.3}{$B.3$})}]
                For any compact set $D \subset C([- h, T], \mathbb{R}^n)$, there exists a number $\lambda > 0$ such that
                \begin{equation*}
                    |H(t, x(\cdot), s) - H(t, y(\cdot), s)|
                    \leq \lambda (1 + \|s\|) \max_{\tau \in [- h, t]} \|x(\tau) - y(\tau)\|
                \end{equation*}
                for every $t \in [0, T]$, $x(\cdot)$, $y(\cdot) \in D$, and $s \in \mathbb{R}^n$.
        \end{description}

        Note that condition \eqref{B.3} implies non-anticipativity of functional \eqref{H_s} for every fixed $s \in \mathbb{R}^n$.
        Note also that, in problem \eqref{HJ}, \eqref{boundary_condition} with Hamiltonian \eqref{Hamiltonian_OCP} associated to optimal control problem \eqref{system}, \eqref{cost_functional}, assumptions \eqref{B.1} to \eqref{B.3} are fulfilled if the functions $f$, $g$, and $\sigma$ satisfy conditions \eqref{A.1} to \eqref{A.3}.

        So, the subject of the paper is problem \eqref{HJ}, \eqref{boundary_condition} under assumptions \eqref{B.1} to \eqref{B.3}.

\section{Minimax Solution}
\label{section_minimax_solution}

    By a classical solution of problem \eqref{HJ}, \eqref{boundary_condition}, it is reasonable to consider a $ci$-smooth functional $\varphi: [0, T] \times C([- h, T], \mathbb{R}^n) \to \mathbb{R}$ that satisfies equation \eqref{HJ} and boundary condition \eqref{boundary_condition}.
    However, under assumptions \eqref{B.1} to \eqref{B.3}, a classical solution of this problem may fail to exist.
    In the theory of Hamilton--Jacobi equations, various notions of generalized solutions are proposed, including minimax (see, e.g., \cite{Subbotin_1991_NA,Subbotin_1995}) and viscosity (see, e.g., \cite{Crandall_Lions_1983,Crandall_Evans_Lions_1984}) solutions.
    This paper continues studies \cite{Krasovskii_Lukoyanov_2000_IMM_Eng,Lukoyanov_2001_DE_Eng,Lukoyanov_2003_1,Lukoyanov_2010_IMM_Eng_1} (see also \cite{Lukoyanov_2011_Eng}) of minimax solutions of problem \eqref{HJ}, \eqref{boundary_condition}.
    Some remarks concerning viscosity solutions of this problem are given in Sect.~\ref{section_viscosity} below.

    \subsection{Definition of Minimax Solution}

        Take the constant $c$ from assumption \eqref{B.2} and, for $(t, x(\cdot)) \in [0, T] \times C([- h, T], \mathbb{R}^n)$, denote
        \begin{equation} \label{Y}
            Y(t, x(\cdot))
            = \big\{ y(\cdot) \in \Lip(t, x(\cdot)):
            \, \|\dot{y}(\tau)\| \leq c (1 + \max_{\xi \in [- h, \tau]} \|y(\xi)\|) \text{ for a.e. } \tau \in [t, T] \big\}.
        \end{equation}

        A functional $\varphi: [0, T] \times C([- h, T], \mathbb{R}^n) \to \mathbb{R}$ is called a {\it minimax solution} of problem \eqref{HJ}, \eqref{boundary_condition} if it is non-anticipative and continuous, satisfies boundary condition \eqref{boundary_condition} and the following property:
        \begin{description}
            \item[\rm{(\mylabel{M}{$M$})}]
                For every $(t, x(\cdot)) \in [0, T) \times C([-h, T], \mathbb{R}^n)$ and $s \in \mathbb{R}^n$, there exists $y(\cdot) \in Y(t, x(\cdot))$ such that
                \begin{equation*}
                    \varphi(\tau, y(\cdot)) - \varphi(t, x(\cdot))
                    = \langle s, y(\tau) - x(t) \rangle - \int_{t}^{\tau} H(\xi, y(\cdot), s) \, \rd \xi,
                    \quad \tau \in [t, T].
                \end{equation*}
        \end{description}

        Property \eqref{M} can be reformulated in terms of (generalized) characteristics of equation \eqref{HJ}.
        Given $(t, x(\cdot)) \in [0, T) \times C([- h, T], \mathbb{R}^n)$ and $s \in \mathbb{R}^n$, consider the Cauchy problem for the functional differential inclusion of retarded type
        \begin{equation} \label{characteristic_DI}
            (\dot{y}(\tau), \dot{z}(\tau))
            \in E(\tau, y(\cdot), s),
            \quad \tau \in [t, T],
        \end{equation}
        where
        \begin{equation} \label{E}
            \begin{array}{c}
                E(\tau, y(\cdot), s)
                = \big\{ (f, g) \in \mathbb{R}^n \times \mathbb{R}: \\[0.5em]
                \displaystyle
                \|f\| \leq c (1 + \max_{\xi \in [- h, \tau]} \|y(\xi)\|),
                \, g = \langle s, f \rangle - H(\tau, y(\cdot), s) \big\},
            \end{array}
        \end{equation}
        and the initial condition
        \begin{equation} \label{initial_condition_DI}
            y(\tau)
            = x(\tau),
            \quad z(\tau)
            = 0,
            \quad \tau \in [- h, t].
        \end{equation}
        By a solution of Cauchy problem \eqref{characteristic_DI}, \eqref{initial_condition_DI}, we mean a pair of functions $(y(\cdot), z(\cdot)) \in C([- h, T], \mathbb{R}^n) \times C([- h, T], \mathbb{R})$ that are Lipschitz continuous on $[t, T]$, meet condition \eqref{initial_condition_DI} and satisfy inclusion \eqref{characteristic_DI} for a.e. $\tau \in [t, T]$.
        Solutions of problem \eqref{characteristic_DI}, \eqref{initial_condition_DI} are called {\it characteristics} of equation \eqref{HJ}.
        The set of all such characteristics is denoted by $CH(t, x(\cdot), s)$.
        Then, property \eqref{M} is equivalent to the following:
        for any $(t, x(\cdot)) \in [0, T) \times C([-h, T], \mathbb{R}^n)$ and $s \in \mathbb{R}^n$, there exists a characteristic $(y(\cdot), z(\cdot)) \in CH(t, x(\cdot), s)$ such that
        \begin{equation} \label{M_equivalent}
            \varphi(\tau, y(\cdot)) - \varphi(t, x(\cdot))
            = z(\tau),
            \quad \tau \in [t, T].
        \end{equation}

        It is often convenient to split the definition of a minimax solution of problem \eqref{HJ}, \eqref{boundary_condition} in two parts and also slightly weaken condition \eqref{M_equivalent}, which leads to the following notions of upper and lower (minimax) solutions of this problem.

        A functional $\varphi: [0, T] \times C([- h, T], \mathbb{R}^n) \to \mathbb{R}$ is called an {\it upper solution} of problem \eqref{HJ}, \eqref{boundary_condition} if it is non-anticipative and lower semicontinuous, satisfies the boundary condition $\varphi(T, x(\cdot)) \geq \sigma(x(\cdot))$, $x(\cdot) \in C([- h, T], \mathbb{R}^n)$, and possesses the property below:
        \begin{description}
            \item[\rm{(\mylabel{M^ast}{$M^\ast$})}]
                For every $(t, x(\cdot)) \in [0, T) \times C([-h, T], \mathbb{R}^n)$, $s \in \mathbb{R}^n$, and $\tau \in (t, T]$, there is a characteristic $(y(\cdot), z(\cdot)) \in CH(t, x(\cdot), s)$ such that
                \begin{equation*}
                    \varphi(\tau, y(\cdot)) - \varphi(t, x(\cdot))
                    \leq z(\tau).
                \end{equation*}
        \end{description}
        Respectively, the functional $\varphi$ is a {\it lower solution} of the problem if it is non-anticipative and upper semicontinuous, meets the condition $\varphi(T, x(\cdot)) \leq \sigma(x(\cdot))$, $x(\cdot) \in C([- h, T], \mathbb{R}^n)$, and has the following property:
        \begin{description}
            \item[\rm{(\mylabel{M_ast}{$M_\ast$})}]
                For every $(t, x(\cdot)) \in [0, T) \times C([-h, T], \mathbb{R}^n)$, $s \in \mathbb{R}^n$, and $\tau \in (t, T]$, there is a characteristic $(y(\cdot), z(\cdot)) \in CH(t, x(\cdot), s)$ such that
                \begin{equation*}
                    \varphi(\tau, y(\cdot)) - \varphi(t, x(\cdot))
                    \geq z(\tau).
                \end{equation*}
        \end{description}

        \begin{proposition}{\rm (see \cite[Proposition~6.1]{Lukoyanov_2003_1})} \label{proposition_upper_lower}
            A functional $\varphi: [0, T] \times C([- h, T], \mathbb{R}^n) \to \mathbb{R}$ is a minimax solution of problem \eqref{HJ}, \eqref{boundary_condition} if and only if it is an upper solution as well as a lower solution of this problem.
        \end{proposition}

    \subsection{Characteristic Complexes}
    \label{subsection_characteristic_complexes}

        In applications, instead of the set-valued function $E$ from \eqref{E}, determining the sets of characteristics of equation \eqref{HJ}, it is often convenient to use another set-valued functions, called characteristic complexes of this equation.

        Let $P$ and $Q$ be non-empty sets.
        Suppose that set-valued functions
        \begin{equation*}
            \begin{array}{c}
                [0, T] \times C([- h, T], \mathbb{R}^n) \times Q \ni (\tau, y(\cdot), q)
                \mapsto E^\ast(\tau, y(\cdot), q) \subset \mathbb{R}^n \times \mathbb{R}, \\[0.5em]
                [0, T] \times C([- h, T], \mathbb{R}^n) \times P \ni (\tau, y(\cdot), p)
                \mapsto E_\ast(\tau, y(\cdot), p) \subset \mathbb{R}^n \times \mathbb{R}
            \end{array}
        \end{equation*}
        satisfy the following conditions:
        \begin{description}
            \item[\rm{(\mylabel{C.1}{$C.1$})}]
                For every $\tau \in [0, T]$, $y(\cdot) \in C([- h, T], \mathbb{R}^n)$, and $q \in Q$, $p \in P$, the sets $E^\ast(\tau, y(\cdot), q)$ and $E_\ast(\tau, y(\cdot), p)$ are non-empty, convex, and compact.

            \item[\rm{(\mylabel{C.2}{$C.2$})}]
                For any $q \in Q$ and $p \in P$, the set-valued functions
                \begin{equation*}
                    \begin{array}{c}
                        [0, T] \times C([- h, T], \mathbb{R}^n) \ni (\tau, y(\cdot))
                        \mapsto E^\ast(\tau, y(\cdot), q) \subset \mathbb{R}^n \times \mathbb{R}, \\[0.5em]
                        [0, T] \times C([- h, T], \mathbb{R}^n) \ni (\tau, y(\cdot))
                        \mapsto E_\ast(\tau, y(\cdot), p) \subset \mathbb{R}^n \times \mathbb{R}
                    \end{array}
                \end{equation*}
                are non-anticipative and upper semicontinuous.

            \item[\rm{(\mylabel{C.3}{$C.3$})}]
                There exists a constant $c > 0$ such that
                \begin{equation*}
                    \max\big\{ \|f\|:
                    \, (f, g) \in E^\ast(\tau, y(\cdot), q) \cup E_\ast(\tau, y(\cdot), p) \big\}
                    \leq c (1 + \max_{\xi \in [- h, \tau]} \|y(\xi)\|)
                \end{equation*}
                for every $\tau \in [0, T]$, $y(\cdot) \in C([- h, T], \mathbb{R}^n)$, and $q \in Q$, $p \in P$.

            \item[\rm{(\mylabel{C.4}{$C.4$})}]
                Given $\tau \in [0, T]$, $y(\cdot) \in C([- h, T], \mathbb{R}^n)$, and $s \in \mathbb{R}^n$, the equalities below hold:
                \begin{equation*}
                    H(\tau, y(\cdot), s)
                    = \sup_{q \in Q} \min_{(f, g) \in E^\ast(\tau, y(\cdot), q)} (\langle s, f \rangle - g)
                    = \inf_{p \in P} \max_{(f, g) \in E_\ast(\tau, y(\cdot), p)} (\langle s, f \rangle - g).
                \end{equation*}
        \end{description}
        Then, the pairs $(Q, E^\ast)$ and $(P, E_\ast)$ are called, respectively, an {\it upper} and a {\it lower characteristic complex} of equation \eqref{HJ}.
        Note that conditions \eqref{C.1} to \eqref{C.3} are rather technical, while condition \eqref{C.4} connects $(Q, E^\ast)$ and $(P, E_\ast)$ with the Hamiltonian $H$.
        Let $\mathcal{E}^\ast(H)$ and $\mathcal{E}_\ast(H)$ denote the collections of all upper and lower characteristic complexes.
        Observe that the set-valued function $E$ from \eqref{E} satisfies the inclusion $(\mathbb{R}^n, E) \in \mathcal{E}^\ast(H) \cap \mathcal{E}_\ast(H)$.

        Let $(t, x(\cdot)) \in [0, T) \times C([- h, T], \mathbb{R}^n)$ and $q \in Q$, $p \in P$ be fixed.
        Consider the Cauchy problem for the functional differential inclusion of retarded type
        \begin{equation*}
            (\dot{y}(\tau), \dot{z}(\tau))
            \in E^\ast(\tau, y(\cdot), q),
            \quad \tau \in [t, T],
        \end{equation*}
        and initial condition \eqref{initial_condition_DI}.
        Denote the set of solutions of this problem by $CH^\ast(t, x(\cdot), q)$.
        Respectively, let $CH_\ast(t, x(\cdot), p)$ be the set of solutions of the Cauchy problem for the differential inclusion
        \begin{equation*}
            (\dot{y}(\tau), \dot{z}(\tau))
            \in E_\ast(\tau, y(\cdot), p),
            \quad \tau \in [t, T],
        \end{equation*}
        and initial condition \eqref{initial_condition_DI}.

        \begin{proposition}{\rm (see \cite[Theorem~5.1]{Lukoyanov_2011_Eng} and also the proof of \cite[Assertion~1]{Lukoyanov_2010_IMM_Eng_1})}
            \label{proposition_characteristic_complexes}
            Let a functional $\varphi: [0, T] \times C([- h, T], \mathbb{R}^n) \to \mathbb{R}$ be non-anticipative and continuous and satisfy boundary condition \eqref{boundary_condition}.
            Then, it is a minimax solution of problem \eqref{HJ}, \eqref{boundary_condition} if and only if the following condition is fulfilled for some $(Q, E^\ast) \in \mathcal{E}^\ast(H)$ and $(P, E_\ast) \in \mathcal{E}_\ast(H)$:
            \begin{description}
                \item[\rm{(\mylabel{MC}{$MC$})}]
                    For every $(t, x(\cdot)) \in [0, T) \times C([- h, T], \mathbb{R}^n)$, $q \in Q$, $p \in P$, and $\tau \in (t, T]$, there exist $(y^\ast(\cdot), z^\ast(\cdot)) \in CH^\ast(t, x(\cdot), q)$ and $(y_\ast(\cdot), z_\ast(\cdot)) \in CH_\ast(t, x(\cdot), p)$ such that
                    \begin{equation*}
                        \varphi(\tau, y^\ast(\cdot)) - \varphi(t, x(\cdot))
                        \leq z^\ast(\tau),
                        \quad \varphi(\tau, y_\ast(\cdot)) - \varphi(t, x(\cdot))
                        \geq z_\ast(\tau).
                    \end{equation*}
            \end{description}
        \end{proposition}

    \subsection{Directional Derivatives of Minimax Solution}

        Property \eqref{M} in the definition of a minimax solution of problem \eqref{HJ}, \eqref{boundary_condition} as well as properties \eqref{M^ast}, \eqref{M_ast}, and \eqref{MC} are non-local, which often complicates their verification.
        In this section, we give an infinitesimal criteria of a minimax solution in terms of appropriate directional derivatives.

        Let $F \subset \mathbb{R}^n$ be a non-empty convex compact set.
        The {\it lower} and {\it upper right derivatives} of a functional $\varphi: [0, T] \times C([- h, T], \mathbb{R}^n) \to \mathbb{R}$ at a point $(t, x(\cdot)) \in [0, T) \times C([- h, T], \mathbb{R}^n)$ in the {\it multi-valued direction $F$} are defined, respectively, by
        \begin{equation*}
            \begin{array}{c}
                \displaystyle
                d^- \{\varphi(t, x(\cdot)) \mid F \}
                = \lim_{\varepsilon \downarrow 0} \inf_{y(\cdot) \in \Omega(t, x(\cdot), F, \varepsilon)} \liminf_{\delta \downarrow 0}
                \big(\varphi(t + \delta, y(\cdot)) - \varphi(t, x(\cdot)) \big) \delta^{- 1}, \\[1em]
                \displaystyle
                d^+ \{\varphi(t, x(\cdot)) \mid F \}
                = \lim_{\varepsilon \downarrow 0} \sup_{y(\cdot) \in \Omega(t, x(\cdot), F, \varepsilon)} \limsup_{\delta \downarrow 0}
                \big(\varphi(t + \delta, y(\cdot)) - \varphi(t, x(\cdot)) \big) \delta^{- 1},
            \end{array}
        \end{equation*}
        where
        \begin{equation*}
            \Omega(t, x(\cdot), F, \varepsilon)
            = \big\{ y(\cdot) \in \Lip(t, x(\cdot)):
            \, \dot{y}(\tau) \in [F]^\varepsilon \text{ for a.e. } \tau \in [t, T] \big\}
        \end{equation*}
        and the symbol $[F]^\varepsilon$ stands for the closed $\varepsilon$-neighborhood of $F$ in $\mathbb{R}^n$.

        \begin{proposition}{\rm (see \cite[Theorem~8.1]{Lukoyanov_2003_1})}
            \label{proposition_directional_derivatives}
            Let a functional $\varphi: [0, T] \times C([- h, T], \mathbb{R}^n) \to \mathbb{R}$ be non-anticipative and continuous and satisfy boundary condition \eqref{boundary_condition}.
            Then, it is a minimax solution of problem \eqref{HJ}, \eqref{boundary_condition} if and only if the differential inequalities below are valid  for every $(t, x(\cdot)) \in [0, T) \times C([- h, T], \mathbb{R}^n)$ and $s \in \mathbb{R}^n$:
            \begin{equation} \label{minimax_directional_-}
                \begin{array}{c}
                    d^- \big\{\varphi(t, x(\cdot)) - \langle s, x(t) \rangle \mid B(t, x(\cdot))\big\} + H(t, x(\cdot), s)
                    \leq 0, \\[0.75em]
                    d^+ \big\{\varphi(t, x(\cdot)) - \langle s, x(t) \rangle \mid B(t, x(\cdot))\big\} + H(t, x(\cdot), s)
                    \geq 0,
                \end{array}
            \end{equation}
            where
            \begin{equation*}
                B(t, x(\cdot))
                = \big\{ f \in \mathbb{R}^n:
                \, \|f\| \leq c (1 + \max_{\tau \in [- h, t]} \|x(\tau)\|) \big\}
            \end{equation*}
            and $c$ is the constant from assumption \eqref{B.2}.
        \end{proposition}

        Note that, if a functional $\varphi: [0, T] \times C([- h, T], \mathbb{R}^n) \to \mathbb{R}$ is $ci$-differentiable at a point $(t, x(\cdot)) \in [0, T) \times C([- h, T], \mathbb{R}^n)$, then (see, e.g., \cite[Proposition~1]{Lukoyanov_2001_PMM_Eng})
        \begin{equation} \label{directional_-_ci-derivatives}
            \begin{array}{c}
                \displaystyle
                d^- \{\varphi(t, x(\cdot)) \mid F \}
                = \partial_t \varphi(t, x(\cdot)) + \min_{f \in F} \langle \nabla \varphi(t, x(\cdot)), f \rangle, \\[1em]
                \displaystyle
                d^+ \{\varphi(t, x(\cdot)) \mid F \}
                = \partial_t \varphi(t, x(\cdot)) + \max_{f \in F} \langle \nabla \varphi(t, x(\cdot)), f \rangle
            \end{array}
        \end{equation}
        for any non-empty convex compact set $F \subset \mathbb{R}^n$.
        More formulas of this kind are given in \cite{Lukoyanov_2001_PMM_Eng,Lukoyanov_2006_IMM_Eng}.
        Note also that, under some additional assumptions on the Hamiltonian $H$ and the boundary functional $\sigma$, a minimax solution of problem \eqref{HJ}, \eqref{boundary_condition} can be characterized with the help of right derivatives in single-valued directions only (see, e.g., \cite{Lukoyanov_2006_IMM_Eng}).

    \subsection{Consistency}

        By virtue of \eqref{directional_-_ci-derivatives}, at any point $(t, x(\cdot)) \in [0, T) \times C([- h, T], \mathbb{R}^n)$ where a functional $\varphi: [0, T] \times C([- h, T], \mathbb{R}^n) \to \mathbb{R}$ is $ci$-differentiable, the pair of inequalities \eqref{minimax_directional_-} is equivalent to equality \eqref{HJ} (see, e.g., \cite[Proposition~12.1]{Lukoyanov_2011_Eng}).
        Thus, Proposition~\ref{proposition_directional_derivatives} allows to conclude that the notion of a minimax solution of problem \eqref{HJ}, \eqref{boundary_condition} is consistent with a classical notion of solution.
        Namely, the following statements are valid.
        \begin{proposition}
            If a continuous functional $\varphi: [0, T] \times C([-h, T], \mathbb{R}^n) \to \mathbb{R}$ is $ci$-differentiable at every point $(t, x(\cdot)) \in [0, T) \times C([- h, T], \mathbb{R}^n)$ and satisfies \eqref{HJ} and \eqref{boundary_condition}, then it is a minimax solution of problem \eqref{HJ}, \eqref{boundary_condition}.
        \end{proposition}
        \begin{proposition}
            If a minimax solution of problem \eqref{HJ}, \eqref{boundary_condition} is $ci$-differentiable at some point $(t, x(\cdot)) \in [0, T) \times C([- h, T], \mathbb{R}^n)$, then it satisfies \eqref{HJ} at this point.
        \end{proposition}

        Let us emphasize that the results given in Sect.~\ref{section_minimax_solution} do not actually require any assumption of kind \eqref{B.3} on Lipschitz continuity of the Hamiltonian.
        They are valid under assumptions \eqref{B.1} and \eqref{B.2} and the supposition that functional \eqref{H_s} is non-anticipative for every $s \in \mathbb{R}^n$.

\section{Well-posedness of Minimax Solution}
\label{section_well_posedeness}

    \subsection{Existence and Uniqueness}
    \label{subsection_Existence_Uniqueness}

        The main contribution of the paper is the following result.
        \begin{theorem} \label{theorem_existence_uniqueness}
            Let assumptions \eqref{B.1} to \eqref{B.3} be satisfied.
            Then, a minimax solution of problem \eqref{HJ}, \eqref{boundary_condition} exists and is unique.
        \end{theorem}

        In accordance with \cite[Theorem~8.1]{Subbotin_1995} (see also \cite{Lukoyanov_2003_1}), the proof of Theorem~\ref{theorem_existence_uniqueness} consists of two parts.
        Firstly, we construct an upper $\varphi^0$ and a lower $\varphi_0$ solutions of problem \eqref{HJ}, \eqref{boundary_condition} such that
        \begin{equation*}
            \varphi^0(t, x(\cdot)) \leq \varphi_0(t, x(\cdot)),
            \quad (t, x(\cdot)) \in [0, T] \times C([- h, T], \mathbb{R}^n).
        \end{equation*}
        Secondly, we establish the so-called {\it comparison principle}, which states that, for any upper $\varphi^\ast$ and any lower $\varphi_\ast$ solutions of the problem, the inequality below holds:
        \begin{equation*}
            \varphi^\ast(t, x(\cdot)) \geq \varphi_\ast(t, x(\cdot)),
            \quad (t, x(\cdot)) \in [0, T] \times C([- h, T], \mathbb{R}^n).
        \end{equation*}
        These two claims yield existence and uniqueness of a minimax solution in view of Proposition~\ref{proposition_upper_lower}.

        Theorem~\ref{theorem_existence_uniqueness} is proved in \cite{Lukoyanov_2003_1} under assumptions \eqref{B.1} and \eqref{B.2} and the following condition, which is stronger than \eqref{B.3}:
        \begin{itemize}
            \item[\rm{(\mylabel{B.4}{$B.4$})}]
                For any compact set $D \subset C([- h, T], \mathbb{R}^n)$, there exists a number $\lambda > 0$ such that
                \begin{equation*}
                    |H(t, x(\cdot), s) - H(t, y(\cdot), s)|
                    \leq \lambda (1 + \|s\|) \bigg( \|x(t) - y(t)\|
                     + \sqrt{\int_{- h}^{t} \|x(\tau) - y(\tau)\|^2 \, \rd \tau} \bigg)
                \end{equation*}
                for every $t \in [0, T]$, $x(\cdot)$, $y(\cdot) \in D$, and $s \in \mathbb{R}^n$.
        \end{itemize}
        In the context of optimal control problems and differential games, condition \eqref{B.4} allows to cover distributed delays only (see \eqref{f_3}).

        Analyzing the proof given in \cite{Lukoyanov_2003_1}, we observe that the first part of the proof does not rely on the specific form of condition \eqref{B.4} and remains valid under assumption \eqref{B.3}.
        Moreover, in the second part of the proof, we can also repeat the reasoning from \cite{Lukoyanov_2003_1} if, given a compact set $D \subset C([- h, T], \mathbb{R}^n)$, we can find a number $\varepsilon_0 > 0$ and a {\it Lyapunov--Krasovskii functional} $\nu_\varepsilon: [0, T] \times C([- h, T], \mathbb{R}^n) \to \mathbb{R}$, depending on the parameter $\varepsilon \in (0, \varepsilon_0]$, such that the statements below are valid:
        \begin{description}
            \item[\rm{(\mylabel{a}{$a$})}]
                The functional $\nu_\varepsilon$ is non-negative and $ci$-smooth.
            \item[\rm{(\mylabel{b}{$b$})}]
                For every $t \in [0, T]$, the estimate $\nu_\varepsilon(t, x(\cdot) \equiv 0) \leq \varepsilon$ is fulfilled.
            \item[\rm{(\mylabel{c}{$c$})}]
                For any constant $C > 0$, the following relation holds:
                \begin{equation*}
                    \max \Big\{ |\sigma(x(\cdot)) - \sigma(y(\cdot))|:
                    \, x(\cdot), y(\cdot) \in D,
                    \, \nu_\varepsilon(T, x(\cdot) - y(\cdot)) \leq C \Big\}
                    \to 0
                    \text{ as } \varepsilon \downarrow 0.
                \end{equation*}
            \item[\rm{(\mylabel{d}{$d$})}]
                The inequality
                \begin{equation*}
                    \partial_t \nu_\varepsilon(t, x(\cdot) - y(\cdot))
                    + H\big(t, x(\cdot), \nabla \nu_\varepsilon(t, x(\cdot) - y(\cdot))\big)
                    - H\big(t, y(\cdot), \nabla \nu_\varepsilon(t, x(\cdot) - y(\cdot))\big)
                    \leq 0
                \end{equation*}
                is valid for every $t \in [0, T)$ and $x(\cdot)$, $y(\cdot) \in D$.
        \end{description}
        The functional $\nu_\varepsilon$ from \cite{Lukoyanov_2003_1} satisfies conditions \eqref{a} to \eqref{d} if \eqref{B.4} is fulfilled.
        In \cite{Lukoyanov_2010_IMM_Eng_1}, an appropriate functional $\nu_\varepsilon$ is constructed under the following assumption, which is weaker than \eqref{B.4} but still stronger than \eqref{B.3}:
        \begin{description}
            \item[\rm{(\mylabel{B.5}{$B.5$})}]
                There exist $J \in \mathbb{N}$ and $\vartheta_j \in (0, h]$, $j \in \overline{1, J}$, such that, for every compact set $D \subset C([- h, T], \mathbb{R}^n)$, there is a number $\lambda > 0$ such that
                \begin{equation*}
                    \begin{array}{c}
                        |H(t, x(\cdot), s) - H(t, y(\cdot), s)|
                        \leq \lambda (1 + \|s\|) \\[0.5em]
                        \displaystyle
                        \times \bigg( \|x(t) - y(t)\|
                        + \sum_{j = 1}^J \|x(t - \vartheta_j) - y(t - \vartheta_j)\|
                        + \sqrt{\int_{- h}^{t} \|x(\tau) - y(\tau)\|^2 \, \rd \tau} \bigg)
                    \end{array}
                \end{equation*}
                for any $t \in [0, T]$, $x(\cdot)$, $y(\cdot) \in D$, and $s \in \mathbb{R}^n$.
        \end{description}
        In applications, assumption \eqref{B.5} allows to consider constant concentrated and distributed delays (see \eqref{f_1} and \eqref{f_3}), however, it does not cover, for example, time-varying delays (see \eqref{f_2}).
        Note that, in general, if we assume only \eqref{B.3}, then the functionals $\nu_\varepsilon$ given in \cite{Lukoyanov_2003_1,Lukoyanov_2010_IMM_Eng_1} may fail to meet condition \eqref{d}.

        Thus, in order to prove Theorem~\ref{theorem_existence_uniqueness} by the scheme from \cite{Lukoyanov_2003_1}, it is sufficient to present a functional $\nu_\varepsilon$ that satisfies \eqref{a} to \eqref{d} under assumptions \eqref{B.1} to \eqref{B.3}.
        Below, we build such a functional based on the constructions from \cite[Sect.~7.5]{Subbotin_1995} and \cite{Zhou_2020_1}.

        Given a compact set $D \subset C([- h, T], \mathbb{R}^n)$, choose a number $\lambda$ by $D$ according to \eqref{B.3}, set $\varkappa = (3 - \sqrt{5}) / 2$, and put $\varepsilon_0 = e^{- \lambda T / \varkappa} / \sqrt{\varkappa}$.
        For every $\varepsilon \in (0, \varepsilon_0]$, define
        \begin{equation*}
            \nu_\varepsilon(t, x(\cdot))
            = \alpha_\varepsilon(t) \beta_\varepsilon(t, x(\cdot)),
            \quad (t, x(\cdot)) \in [0, T] \times C([- h, T], \mathbb{R}^n),
        \end{equation*}
        where
        \begin{equation*}
            \alpha_\varepsilon(t)
            = (e^{- \lambda t / \varkappa} - \varepsilon \sqrt{\varkappa}) \varepsilon^{- 1},
            \quad \beta_\varepsilon(t, x(\cdot))
            = \sqrt{\varepsilon^4 + V(t, x(\cdot))},
        \end{equation*}
        and
        \begin{equation} \label{V}
            V(t, x(\cdot))
            = \begin{cases}
                \displaystyle
                \frac{\big(\|x(\cdot)\|_{[- h, t]}^2 - \|x(t)\|^2\big)^2}{\|x(\cdot)\|_{[- h, t]}^2} + \|x(t)\|^2
                & \mbox{if } \|x(\cdot)\|_{[- h, t]} > 0, \\
                0 & \mbox{if } \|x(\cdot)\|_{[- h, t]} = 0.
              \end{cases}
        \end{equation}
        Here, in accordance with \eqref{norm_C} and \eqref{x_t}, we denote
        \begin{equation*}
            \|x(\cdot)\|_{[- h, t]}
            = \|x_t(\cdot)\|_{[- h, t]}
            = \max_{\tau \in [- h, t]} \|x(\tau)\|.
        \end{equation*}

        For every $\varepsilon \in (0, \varepsilon_0]$, the functional $\nu_\varepsilon$ is non-negative by the choice of $\varepsilon_0$.
        Further, it can be verified directly (see \cite[Lemma~2.3]{Zhou_2020_1} and also Sect.~\ref{subsection_derivatives} below) that the functional $V$ from \eqref{V} is $ci$-smooth and its $ci$-derivatives are given by
        \begin{equation*}
            \partial_t V(t, x(\cdot))
            = 0,
            \ \ \nabla V(t, x(\cdot))
            = \begin{cases}
                \displaystyle
                \bigg(2 - \frac{4 \big(\|x(\cdot)\|_{[- h, t]}^2 - \|x(t)\|^2\big)}{\|x(\cdot)\|_{[- h, t]}^2} \bigg) x(t)
                & \hspace*{-0.3em} \mbox{if } \|x(\cdot)\|_{[- h, t]} > 0, \\
                0 & \hspace*{-0.3em} \mbox{if } \|x(\cdot)\|_{[- h, t]} = 0
            \end{cases}
        \end{equation*}
        for every $(t, x(\cdot)) \in [0, T) \times C([- h, T], \mathbb{R}^n)$.
        Hence, the functional $\nu_\varepsilon$ is $ci$-smooth and
        \begin{equation} \label{nu_varepsilon_derivatives}
            \partial_t \nu_\varepsilon(t, x(\cdot))
            = - \frac{\lambda e^{- \lambda t / \varkappa} \beta_\varepsilon(t, x(\cdot))}{\varkappa \varepsilon},
            \quad \nabla \nu_\varepsilon(t, x(\cdot))
            = \frac{\alpha_\varepsilon(t)}{2 \beta_\varepsilon(t, x(\cdot))} \nabla V(t, x(\cdot)).
        \end{equation}
        So, condition \eqref{a} is satisfied.
        The validity of \eqref{b} follows from the equality $V(t, x(\cdot) \equiv 0) = 0$, $t \in [0, T]$.
        Condition \eqref{c} is fulfilled due to continuity of the boundary functional $\sigma$ (see assumption \eqref{B.1}) and the inequality (see \cite[Lemma~2.3]{Zhou_2020_1})
        \begin{equation} \label{V_bound}
            V(t, x(\cdot))
            \geq \varkappa \|x(\cdot)\|_{[- h, t]}^2,
            \quad (t, x(\cdot)) \in [0, T] \times C([- h, T], \mathbb{R}^n).
        \end{equation}
        The last condition \eqref{d} holds by virtue of the choice of $\lambda$, relations \eqref{nu_varepsilon_derivatives} and \eqref{V_bound}, and the estimate
        \begin{equation*}
            \|\nabla V(t, x(\cdot))\|
            \leq 2 \|x(t)\|,
            \quad (t, x(\cdot)) \in [0, T) \times C([- h, T], \mathbb{R}^n).
        \end{equation*}
        Thus, the proposed functional $\nu_\varepsilon$ meets all conditions \eqref{a} to \eqref{d} under assumptions \eqref{B.1} to \eqref{B.3}, and, consequently, we conclude that Theorem~\ref{theorem_existence_uniqueness} is proved.

        The above reasoning show that, instead of condition \eqref{B.3}, it can only be assumed that there exists a functional $\nu_\varepsilon$ satisfying \eqref{a} to \eqref{d} (see also \cite{Lukoyanov_2006_IMM_Eng}).
        Moreover, according to \cite{Lukoyanov_2004_PMM_Eng}, the proposed functional $\nu_\varepsilon$ can also be used as the basis for constructing optimal positional (closed-loop) strategies in zero-sum differential games for functional differential equations of retarded type.

        Let us mention that, in \cite{Bayraktar_Keller_2018}, a Lipschitz continuity condition of a slightly another kind is considered for the Hamiltonian $H$.
        In relation to problem \eqref{HJ}, \eqref{boundary_condition}, this condition can be formulated as follows:
        \begin{description}
            \item[\rm{(\mylabel{B.6}{$B.6$})}]
                For every $L \geq 0$ and $(t, x(\cdot)) \in [0, T) \times C([- h, T], \mathbb{R}^n)$, there is a modulus $m$ such that
                \begin{equation*}
                    \begin{array}{c}
                        \big| H\big(\tau, y(\cdot), (y(\tau) - z(\tau)) \varepsilon^{- 1} \big)
                        - H \big(\tau, z(\cdot), (y(\tau) - z(\tau)) \varepsilon^{- 1} \big) \big| \\[0.5em]
                        \displaystyle
                        \leq m \Big( \|y(\tau) - z(\tau)\|^2 \varepsilon^{- 1} + \max_{\xi \in [- h, \tau]} \|y(\xi) - z(\xi)\| \Big)
                    \end{array}
                \end{equation*}
                for any $\varepsilon > 0$, $\tau \in [t, T]$, and $y(\cdot)$, $z(\cdot) \in Y_L(t, x(\cdot))$.
                Here, the set $Y_L(t, x(\cdot))$ is defined by \eqref{Y} with the constant $c$ replaced by $L$.
        \end{description}
        Note that, for example, the Hamiltonians
        \begin{equation*}
            H(t, x(\cdot), s)
            = \langle s, x(t - h) \rangle,
            \quad H(t, x(\cdot), s)
            = \langle s, x(t / 2) \rangle,
            \quad H(t, x(\cdot), s)
            = \Big\langle s, \int_{- h}^{t} x(\tau) \, \rd \tau \Big\rangle,
        \end{equation*}
        which obviously satisfy \eqref{B.3}, do not satisfy \eqref{B.6}.
        So, condition \eqref{B.3} does not imply condition \eqref{B.6}.

    \subsection{Stability}

        Let, for every $k \in \mathbb{N} \cup \{0\}$, a Hamiltonian $H_k: [0, T] \times C([- h, T], \mathbb{R}^n) \times \mathbb{R}^n \to \mathbb{R}$ and a boundary functional $\sigma_k: C([- h, T], \mathbb{R}^n) \to \mathbb{R}$ be given satisfying assumptions \eqref{B.1}, \eqref{B.2} with the constant $c$ independent on $k$, and \eqref{B.3}.
        Suppose that, for any compact set $D \subset C([- h, T], \mathbb{R}^n)$ and any $s \in \mathbb{R}^n$, we have
        \begin{equation*}
            |H_0(t, x(\cdot), s) - H_k(t, x(\cdot), s)| + |\sigma_0(x(\cdot)) - \sigma_k(x(\cdot))| \to 0
            \text{ as } k \to \infty
        \end{equation*}
        uniformly in $t \in [0, T]$ and $x(\cdot) \in D$.

        For every $k \in \mathbb{N} \cup \{0\}$, consider the minimax solution $\varphi_k: [0, T] \times C([- h, T], \mathbb{R}^n) \to \mathbb{R}$ of the Cauchy problem for the Hamilton--Jacobi equation with $ci$-derivatives
        \begin{equation*}
            \partial_t \varphi(t, x(\cdot)) + H_k \big( t, x(\cdot), \nabla \varphi(t, x(\cdot)) \big)
            = 0,
            \quad (t, x(\cdot)) \in [0, T) \times C([- h, T], \mathbb{R}^n),
        \end{equation*}
        and the boundary condition
        \begin{equation*}
            \varphi(T, x(\cdot))
            = \sigma_k(x(\cdot)),
            \quad x(\cdot) \in C([- h, T], \mathbb{R}^n).
        \end{equation*}

        Owing to the results of Sect.~\ref{subsection_Existence_Uniqueness}, the following statement can be proved by the scheme from \cite{Krasovskii_Lukoyanov_2000_IMM_Eng,Lukoyanov_2001_DE_Eng} (see also \cite[Theorem~9.1]{Lukoyanov_2011_Eng}).
        \begin{proposition}
            Under the suppositions made above, it holds that $\varphi_k(t, x(\cdot)) \to \varphi_0(t, x(\cdot))$ as $k \to \infty$ uniformly in $t \in [0, T]$ and $x(\cdot) \in D$ for every compact set $D \subset C([- h, T], \mathbb{R}^n)$.
        \end{proposition}

    \subsection{Consistency with Non-Path-Dependent Case}

        Suppose that the Hamiltonian $H$ and the boundary functional $\sigma$ in problem \eqref{HJ}, \eqref{boundary_condition} are of the form
        \begin{equation} \label{overline_H_sigma}
            \begin{array}{c}
                H(t, x(\cdot), s)
                = \widehat{H}(t, x(t), s),
                \quad \sigma(x(\cdot))
                = \widehat{\sigma}(x(T)), \\[0.5em]
                t \in [0, T], \quad x(\cdot) \in C([- h, T], \mathbb{R}^n), \quad s \in \mathbb{R}^n,
            \end{array}
        \end{equation}
        for some functions $\widehat{H}: [0, T] \times \mathbb{R}^n \times \mathbb{R}^n \to \mathbb{R}$ and $\widehat{\sigma}: \mathbb{R}^n \to \mathbb{R}$.
        Consider the Cauchy problem for the Hamilton--Jacobi equation with partial derivatives
        \begin{equation} \label{HJ_ordinary}
            \frac{\partial \widehat{\varphi}}{\partial t} (t, x) + \widehat{H}(t, x, \nabla_x \widehat{\varphi}(t, x))
            = 0,
            \quad (t, x) \in (0, T) \times \mathbb{R}^n,
        \end{equation}
        and the boundary condition
        \begin{equation} \label{boundary_condition_ordinary}
            \widehat{\varphi}(T, x)
            = \widehat{\sigma}(x),
            \quad x \in \mathbb{R}^n.
        \end{equation}
        Conditions \eqref{B.1} to \eqref{B.3} imply that the functions $\widehat{H}$ and $\widehat{\sigma}$ possess the properties providing that a minimax solution $\widehat{\varphi}: [0, T] \times \mathbb{R}^n \to \mathbb{R}$ of problem \eqref{HJ_ordinary}, \eqref{boundary_condition_ordinary} exists and is unique (see, e.g., \cite[Theorem~8.1]{Subbotin_1995}).
        In particular, the function $\widehat{\varphi}$ is continuous and satisfies \eqref{boundary_condition_ordinary} and the following property (see, e.g., \cite[Definition~(M.2) in Sect.~6.3]{Subbotin_1995}):
        \begin{description}
            \item[\rm{(\mylabel{widehat_M}{$\widehat{M}$})}]
                For every $(t, x) \in [0, T) \times \mathbb{R}^n$ and $s \in \mathbb{R}^n$, there exists a Lipschitz continuous function $y: [t, T] \to \mathbb{R}^n$ such that
                \begin{equation*}
                    y(t)
                    = x,
                    \quad \|\dot{y}(\tau)\|
                    \leq c (1 + \|y(\tau)\|)
                    \text{ for a.e. } \tau \in [t, T],
                \end{equation*}
                and the equality below holds:
                \begin{equation*}
                    \widehat{\varphi}(\tau, y(\tau)) - \widehat{\varphi}(t, x)
                    = \langle s, y(\tau) - x \rangle - \int_{t}^{\tau} \widehat{H}(\xi, y(\xi), s) \, \rd \xi,
                    \quad \tau \in [t, T].
                \end{equation*}
        \end{description}
        Here, the constant $c$ can be taken from assumption \eqref{B.2}.
        Then, we obtain that the functional
        \begin{equation} \label{overline_varphi}
            \varphi(t, x(\cdot))
            = \widehat{\varphi}(t, x(t)),
            \quad (t, x(\cdot)) \in [0, T] \times C([- h, T], \mathbb{R}^n),
        \end{equation}
        is non-anticipative and continuous, satisfies \eqref{boundary_condition} and property \eqref{M}.
        So, in view of Theorem~\ref{theorem_existence_uniqueness}, the following result is valid.
        \begin{proposition}
            Let the Hamiltonian $H$ and the boundary functional $\sigma$ in problem \eqref{HJ}, \eqref{boundary_condition} be of form \eqref{overline_H_sigma}.
            Then, the minimax solution of this problem coincides with the functional $\varphi$ from \eqref{overline_varphi}.
        \end{proposition}

\section{Remarks on Problem Statement}
\label{section_remarks_on_problem_statement}

    \subsection{Continuous Non-Anticipative Functionals}
    \label{subsection_continuity}

        In this paper, following, e.g., \cite{Lukoyanov_2007_IMM_Eng,Lukoyanov_2010_IMM_Eng_1}, we consider problem \eqref{HJ}, \eqref{boundary_condition} in the set $[0, T] \times C([- h, T], \mathbb{R}^n)$ assuming that the sought functional $\varphi: [0, T] \times C([- h, T], \mathbb{R}^n) \to \mathbb{R}$ is non-anticipative and at least continuous with respect to the metric $\dist$ from \eqref{dist}.

        Another widely accepted approach (see, e.g., \cite{Pham_Zhang_2014,Kaise_2015,Tang_Zhang_2015,Ekren_Touzi_Zhang_2016_1,Bayraktar_Keller_2018,Kaise_Kato_Takahashi_2018,Saporito_2019,Cosso_Russo_2019,Zhou_2020_1,Zhou_2020_2}) is to endow the set $[0, T] \times C([- h, T], \mathbb{R}^n)$ with the pseudometric
        \begin{equation} \label{dist_infty}
            \rho\big((t, x(\cdot)), (\tau, y(\cdot)) \big)
            = |t - \tau| + \|x(\cdot \wedge t) - y(\cdot \wedge \tau)\|_{[- h, T]},
        \end{equation}
        where $a \wedge b = \min\{a, b\}$ for $a$, $b \in \mathbb{R}$ and, respectively,
        \begin{equation*}
            x(\xi \wedge t)
            = \begin{cases}
                x(\xi) & \mbox{if } \xi \in [- h, t), \\
                x(t) & \mbox{if } \xi \in [t, T].
              \end{cases}
        \end{equation*}
        Let us observe that a functional $\varphi: [0, T] \times C([- h, T], \mathbb{R}^n) \to \mathbb{R}$ is continuous with respect to the pseudometric $\rho$ if and only if it is non-anticipative and continuous with respect to the metric $\dist$.
        So, the results of this paper remain valid if we reformulate problem \eqref{HJ}, \eqref{boundary_condition} in the set $[0, T] \times C([- h, T], \mathbb{R}^n)$ with the pseudometric $\rho$ without mentioning any non-anticipativity properties.

        Further, based on a standard procedure for passing from a pseudometric space to the induced metric space of equivalence classes (see, e.g., \cite[Ch.~4, Theorem~15]{Kelley_1975}), we can move from the set $[0, T] \times C([- h, T], \mathbb{R}^n)$ with the pseudometric $\rho$ to the set
        \begin{equation*}
            G
            = \big\{(t, w(\cdot)):
            \, t \in [0, T], \, w(\cdot) \in C([- h, t], \mathbb{R}^n) \big\}
        \end{equation*}
        with the metric
        \begin{equation*}
            \rho_\ast\big( (t, w(\cdot)), (\tau, r(\cdot)) \big)
            = |t - \tau| + \|w(\cdot \wedge t) - r(\cdot \wedge \tau)\|_{[- h ,T]},
            \quad (t, w(\cdot)), (\tau, r(\cdot)) \in G.
        \end{equation*}
        In this way, we can identify every non-anticipative functional $\varphi: [0, T] \times C([- h, T], \mathbb{R}^n) \to \mathbb{R}$ with the functional $\psi: G \to \mathbb{R}$ given by
        \begin{equation*}
            \psi(t, w(\cdot))
            = \varphi(t, w(\cdot \wedge t)),
            \quad (t, w(\cdot)) \in G.
        \end{equation*}
        Then, the functional $\varphi$ is continuous with respect to the pseudometric $\rho$ (or, which is the same, it is non-anticipative and continuous with respect to the metric $\dist$) if and only if the functional $\psi$ is continuous with respect to the metric $\rho_\ast$.
        Thus, problem \eqref{HJ}, \eqref{boundary_condition} can also be reformulated in the set $G$ with the metric $\rho_\ast$, when the Hamiltonian $H$ is defined on triples $(t, w(\cdot), s)$ and the unknown is a functional from $G$ to $\mathbb{R}$.

        Let us note that, in, e.g., \cite{Krasovskii_Lukoyanov_2000_IMM_Eng,Lukoyanov_2001_DE_Eng,Lukoyanov_2003_1,Lukoyanov_2011_Eng}, problem \eqref{HJ}, \eqref{boundary_condition} is considered in the set $G$ but with the metric
        \begin{equation*}
            \rho_H\big( (t, w(\cdot)), (\tau, r(\cdot)) \big)
            = \max \Big\{ \rho_H^\ast\big( (t, w(\cdot)), (\tau, r(\cdot)) \big),
            \rho_H^\ast \big( (\tau, r(\cdot)), (t, w(\cdot)) \big) \Big\},
        \end{equation*}
        where
        \begin{equation*}
            \rho_H^\ast \big( (t, w(\cdot)), (\tau, r(\cdot)) \big)
            = \max_{\xi \in [0, t]} \min_{\eta \in [0, \tau]}
            \sqrt{|\xi - \eta|^2 + \|w(\xi) - r(\eta)\|^2}.
        \end{equation*}
        Actually, $\rho_H((t, w(\cdot)), (\tau, r(\cdot)))$ is the Hausdorff distance between the graphics of the functions $w(\cdot)$ and $r(\cdot)$ as compact subsets of $\mathbb{R}^{n + 1}$.
        In accordance with, e.g., \cite[Ch.~2]{Sendov_1990}, the metrics $\rho_H$ and $\rho_\ast$ are not strongly equivalent, but they induce the same topology on $G$, which means that, if $(t_k, w_k(\cdot)) \in G$, $k \in \mathbb{N} \cup \{0\}$, then $\rho_H ((t_0, w_0(\cdot)), (t_k, w_k(\cdot))) \to 0$ as $k \to \infty$ if and only if $\rho_\ast((t_0, w_0(\cdot)), (t_k, w_k(\cdot))) \to 0$ as $k \to \infty$.
        Respectively, any functional $\psi: G \to \mathbb{R}$ is continuous with respect to the metric $\rho_H$ if and only if it is continuous with respect to the metric $\rho_\ast$.

        To conclude this section, let us note that the arguments above allow us to apply in this paper the results from \cite{Krasovskii_Lukoyanov_2000_IMM_Eng,Lukoyanov_2001_DE_Eng,Lukoyanov_2001_PMM_Eng,Lukoyanov_2003_1,Lukoyanov_2004_PMM_Eng,Lukoyanov_2006_IMM_Eng,Lukoyanov_2011_Eng}.

    \subsection{Derivatives of Non-Anticipative Functionals}
    \label{subsection_derivatives}

        In this paper, we consider path-dependent Hamilton--Jacobi equation \eqref{HJ} that involves co-invariant ($ci$-) derivatives, defined by relation \eqref{ci_derivatives_definition}.

        Invariant and co-invariant derivatives of non-anticipative functionals were introduced in \cite{Kim_1985} for studying stability issues of functional differential equations on the basis of Lyapunov--Krasovskii functionals (see, e.g., \cite{Krasovskii_1963}).
        Later, the corresponding technique was developed and applied in various parts of the theory of functional differential equations (see, e.g., \cite{Kim_1999} and the references therein).
        In \cite{Kim_1999}, functionals of the following form are investigated:
        \begin{equation*}
            [0, T] \times \mathbb{R}^n \times PC([- h, T], \mathbb{R}^n) \ni (t, z, x(\cdot)) \mapsto \Phi(t, z, x(\cdot)) \in \mathbb{R},
        \end{equation*}
        where $PC([- h, T], \mathbb{R}^n)$ is the set of piecewise continuous functions $x: [-h, T] \to \mathbb{R}^n$.
        For such a functional $\Phi$, the {\it $ci$-derivative} $\partial \Phi(t, z, x(\cdot)) \in \mathbb{R}$ at a point $(t, z, x(\cdot)) \in [0, T) \times \mathbb{R}^n \times PC([- h, T], \mathbb{R}^n)$ is defined as
        \begin{equation*}
            \partial \Phi(t, z, x(\cdot))
            = \lim_{\delta \downarrow 0} \big( \Phi(t + \delta, z, y(\cdot)) - \Phi(t, z, x(\cdot)) \big) \delta^{- 1}
        \end{equation*}
        if the limit exists and does not depend on $y(\cdot) \in \Lip(t, z, x(\cdot))$, where $\Lip(t, z, x(\cdot))$ is the set of $y(\cdot) \in PC([- h, T], \mathbb{R}^n)$ such that $y(\tau) = x(\tau)$, $\tau \in [- h, t)$, $y(t) = z$, and $y(\cdot)$ is Lipschitz continuous on $[t, T]$.
        Concerning the variable $z$, the usual gradient $\nabla_z \Phi(t, z, x(\cdot)) \in \mathbb{R}^n$ is considered.
        Note that this approach requires the functional $\Phi$ to be defined on piecewise continuous functions $x(\cdot)$, and, in contrast to \eqref{ci_derivatives_definition}, it can not be applied directly for functionals defined on continuous functions $x(\cdot)$ only.
        Nevertheless, if, given a functional $\Phi$, we consider the functional
        \begin{equation*}
            \varphi(t, x(\cdot))
            = \Phi(t, x(t), x(\cdot)),
            \quad (t, x(\cdot)) \in [0, T] \times C([- h, T], \mathbb{R}^n),
        \end{equation*}
        then, under certain assumptions, we obtain
        \begin{equation*}
            \partial_t \varphi(t, x(\cdot))
            = \partial \Phi(t, x(t), x(\cdot)),
            \quad \nabla \varphi(t, x(\cdot))
            = \nabla_z \Phi(t, x(t), x(\cdot)).
        \end{equation*}
        In view of this connection, for the values $\partial_t \varphi(t, x(\cdot))$ and $\nabla \varphi(t, x(\cdot))$ given by \eqref{ci_derivatives_definition}, we also use the term $ci$-derivatives.

        Let us mention that a close approach to differentiation of non-anticipative functionals was proposed in \cite{Aubin_Haddad_2002}, where the corresponding derivatives are called Clio derivatives.

        Path-dependent Hamilton--Jacobi equations are often considered (see, e.g., \cite{Tang_Zhang_2015,Saporito_2019,Cosso_Russo_2019,Zhou_2020_1,Zhou_2020_2}) with horizontal and vertical derivatives \cite{Dupire_2009}.
        Let $D([- h, T], \mathbb{R}^n)$ be the Banach space of c\`{a}dl\`{a}g functions $x: [- h, T] \to \mathbb{R}$ with the norm
        \begin{equation*}
            \|x(\cdot)\|_{[- h, T]}
            = \sup_{t \in [- h, T]} \|x(t)\|,
            \quad x(\cdot) \in D([- h, T], \mathbb{R}^n).
        \end{equation*}
        Consider the set $[0, T] \times D([- h, T], \mathbb{R}^n)$ endowed with the pseudometric $\rho$ given by \eqref{dist_infty}.
        The {\it horizontal derivative} $\partial^H \Phi(t, x(\cdot)) \in \mathbb{R}$ of a functional $\Phi: [0, T] \times D([- h, T], \mathbb{R}^n) \to \mathbb{R}$ at a point $(t, x(\cdot)) \in [0, T) \times D([- h, T], \mathbb{R}^n)$ is defined as
        \begin{equation*}
            \partial^H \Phi(t, x(\cdot))
            = \lim_{\delta \downarrow 0} \big( \Phi(t + \delta, x(\cdot \wedge t)) - \Phi(t, x(\cdot)) \big) \delta^{- 1}.
        \end{equation*}
        The {\it vertical derivative} $\partial^V \Phi(t, x(\cdot)) \in \mathbb{R}^n$ of $\Phi$ at a point $(t, x(\cdot)) \in [0, T] \times D([- h, T], \mathbb{R}^n)$ is defined as $\partial^V \Phi(t, x(\cdot)) = (\partial_1^V \Phi(t, x(\cdot)), \ldots, \partial_n^V \Phi(t, x(\cdot))) \in \mathbb{R}^n$, where
        \begin{equation*}
            \partial_i^V \Phi(t, x(\cdot))
            = \lim_{\delta \to 0} \big( \Phi(t, x(\cdot) + \delta e_i 1_{[t, T]}(\cdot)) - \Phi(t, x(\cdot)) \big) \delta^{- 1}.
        \end{equation*}
        Here, $e_i \in \mathbb{R}^n$, $i \in \overline{1, n}$, is the standard orthonormal basis of $\mathbb{R}^n$, and $1_{[t, T]}(\cdot) \in D([- h, T], \mathbb{R})$ stands for the indicator function of the segment $[t, T]$.
        Note again that the vertical derivatives can not be defined directly for functionals $\varphi: [0, T] \times C([- h, T], \mathbb{R}^n) \to \mathbb{R}$ (in this regard, see, e.g., \cite{Cosso_Russo_2019}).
        The following statement is an analogue of Proposition~\ref{proposition_ci_smooth}.
        \begin{proposition}{\rm (see \cite[Theorem~2.1]{Cosso_Russo_2019} and also \cite[Theorem~2.6]{Zhou_2020_1})}
        \label{proposition_HV_derivatives}
            Suppose that a functional $\Phi: [0, T] \times D([- h, T], \mathbb{R}^n) \to \mathbb{R}$ is continuous and possesses a horizontal derivative at every point $(t, x(\cdot)) \in [0, T) \times D([- h, T], \mathbb{R}^n)$ and a vertical derivative at every point $(t, x(\cdot)) \in [0, T] \times D([- h, T], \mathbb{R}^n)$, and let these derivatives
            \begin{equation*}
                \begin{array}{c}
                    [0, T) \times D([- h, T], \mathbb{R}^n) \ni (t, x(\cdot)) \mapsto \partial^H \Phi(t, x(\cdot)) \in \mathbb{R}, \\[0.5em]
                    [0, T] \times D([- h, T], \mathbb{R}^n) \ni (t, x(\cdot)) \mapsto \partial^V \Phi(t, x(\cdot)) \in \mathbb{R}^n
                \end{array}
            \end{equation*}
            be continuous.
            Then, for every $(t, x(\cdot)) \in [0, T) \times C([- h, T], \mathbb{R}^n)$ and $y(\cdot) \in \Lip(t, x(\cdot))$, the equality below holds for any $\tau \in [t, T)$:
            \begin{equation} \label{proposition_functional_derivatives_main}
                \Phi(\tau, y(\cdot))
                = \Phi(t, x(\cdot))
                + \int_{t}^{\tau} \partial^H \Phi(\xi, y(\cdot)) \, \rd \xi
                + \int_{t}^{\tau} \langle \partial^V \Phi(\xi, y(\cdot)), \dot{y}(\xi) \rangle \, \rd \xi.
            \end{equation}
        \end{proposition}

        As a consequence, we obtain that, if a functional $\Phi: [0, T] \times D([- h, T], \mathbb{R}^n) \to \mathbb{R}$ satisfies the assumptions of Proposition~\ref{proposition_HV_derivatives}, then its restriction $\varphi$ to $[0, T] \times C([- h, T], \mathbb{R}^n)$ is a $ci$-smooth functional and, for every $(t, x(\cdot)) \in [0, T) \times C([- h, T], \mathbb{R}^n)$,
        \begin{equation*}
            \partial_t \varphi(t, x(\cdot))
            = \partial^H \Phi(t, x(\cdot)),
            \quad \nabla \varphi(t, x(\cdot))
            = \partial^V \Phi(t, x(\cdot)).
        \end{equation*}
        Taking this connection into account, we can apply the results of the paper to study classical and minimax solutions of problem \eqref{HJ}, \eqref{boundary_condition} reformulated in terms of horizontal and vertical derivatives.

        Let us mention that there is a relation between horizontal and vertical derivatives and Fr\'{e}chet derivatives of non-anticipative functionals (see, e.g., \cite{Ji_Yang_2015}).
        Nevertheless, to the best of our knowledge, it remains unclear how to apply the theory of Hamilton--Jacobi equations with Fr\'{e}chet derivatives in infinite dimensional spaces to investigate problem \eqref{HJ}, \eqref{boundary_condition} under assumptions \eqref{B.1} to \eqref{B.3}.

        Finally, we observe that, introducing derivatives of non-anticipative functionals in one way or another, it is important to have a formula of kind \eqref{proposition_ci_smooth_main} or \eqref{proposition_functional_derivatives_main}.
        Therefore, sometimes (see, e.g., \cite{Pham_Zhang_2014,Ekren_Touzi_Zhang_2016_1,Bayraktar_Keller_2018}), such formulas are taken directly as the basis for definition of the corresponding derivatives.

\section{Remarks on the Case of Homogenous Hamiltonian}
\label{section_remarks_homogenous}

    In, e.g., \cite{Krasovskii_Lukoyanov_2000_IMM_Eng,Lukoyanov_2001_DE_Eng}, problem \eqref{HJ}, \eqref{boundary_condition} is studied under Lipschitz continuity assumptions of kind \eqref{B.3} and some additional suppositions concerning positive homogeneity of the Hamiltonian $H$ with respect to the impulse variable $s$.

    Namely, in \cite{Lukoyanov_2001_DE_Eng}, problem \eqref{HJ}, \eqref{boundary_condition} is considered under assumptions \eqref{B.1}, \eqref{B.2}, and
    \begin{description}
        \item[\rm{(\mylabel{B.7}{$B.7$})}]
            For every compact set $D \subset C([- h, T], \mathbb{R}^n)$, there exists a number $\lambda > 0$ such that
            \begin{equation*}
                |H(t, x(\cdot), s) - H(t, y(\cdot), s)|
                \leq \lambda \max_{\tau \in [- h, t]} \|x(\tau) - y(\tau)\|
            \end{equation*}
            for any $t \in [0, T]$, $x(\cdot)$, $y(\cdot) \in D$, and $s \in \mathbb{R}^n$ such that $\|s\| = 1$.

        \item[\rm{(\mylabel{B.8}{$B.8$})}]
            For every $t \in [0, T]$, $x(\cdot) \in C([0, T], \mathbb{R}^n)$, $s \in \mathbb{R}^n$, and $\alpha \geq 0$, the equality below holds:
            \begin{equation*}
                H(t, x(\cdot), \alpha s)
                = \alpha H(t, x(\cdot), s).
            \end{equation*}
    \end{description}

    Note that assumptions \eqref{B.7} and \eqref{B.8} imply condition \eqref{B.3}, while conditions \eqref{B.4} and \eqref{B.5} may fail to be valid.
    Note also that, for Hamiltonian \eqref{Hamiltonian_OCP} associated to optimal control problem \eqref{system}, \eqref{cost_functional}, assumption \eqref{B.8} is fulfilled if $g(\tau, y(\cdot), u) \equiv 0$.

    Given $(\tau, y(\cdot)) \in [0, T] \times C([- h, T], \mathbb{R}^n)$ and $s \in \mathbb{R}^n$, denote
    \begin{equation*}
        \begin{array}{c}
            E^\ast(\tau, y(\cdot), s)
            = \big\{ f \in F(\tau, y(\cdot)):
            \, \langle s, f \rangle \geq H(\tau, y(\cdot), s) \big\} \times \{0\} \subset \mathbb{R}^n \times \mathbb{R}, \\[0.5em]
            E_\ast(\tau, y(\cdot), s)
            = \big\{ f \in F(\tau, y(\cdot)):
            \, \langle s, f \rangle \leq H(\tau, y(\cdot), s) \big\} \times \{0\}  \subset \mathbb{R}^n \times \mathbb{R},
        \end{array}
    \end{equation*}
    where
    \begin{equation*}
        F(\tau, y(\cdot))
        = \big\{ f \in \mathbb{R}^n:
        \, \|f\| \leq \sqrt{2} c (1 + \max_{\xi \in [- h, \tau]} \|y(\xi)\|) \big\}
    \end{equation*}
    and $c$ is the constant from assumption \eqref{B.2}.
    Observe that the set-valued functions $E^\ast$ and $E_\ast$ satisfy the inclusions $(\mathbb{R}^n, E^\ast) \in \mathcal{E}^\ast(H)$ and $(\mathbb{R}^n, E_\ast) \in \mathcal{E}_\ast(H)$ (see Sect.~\ref{subsection_characteristic_complexes}).

    By \cite[Theorem~5.1]{Lukoyanov_2001_DE_Eng}, there exists a unique functional $\varphi: [0, T] \times C([- h, T], \mathbb{R}^n) \to \mathbb{R}$ that is non-anticipative and continuous, meets boundary condition \eqref{boundary_condition} and possesses property \eqref{MC} for $(\mathbb{R}^n, E^\ast)$ and $(\mathbb{R}^n, E_\ast)$.
    Then, by Proposition~\ref{proposition_characteristic_complexes} and Theorem~\ref{theorem_existence_uniqueness}, this functional $\varphi$ coincides with the minimax solution of problem \eqref{HJ}, \eqref{boundary_condition}.
    Thus, the present paper generalizes the corresponding results obtained for problem \eqref{HJ}, \eqref{boundary_condition} under assumptions \eqref{B.1}, \eqref{B.2}, \eqref{B.7}, and \eqref{B.8}.

    In \cite{Krasovskii_Lukoyanov_2000_IMM_Eng}, the case is investigated when the Hamiltonian $H$ is not homogeneous, but problem \eqref{HJ}, \eqref{boundary_condition} can be reduced to an auxiliary problem with a homogenous Hamiltonian $\overline{H}$.
    It is supposed that assumption \eqref{B.1} is fulfilled and the following conditions hold:
    \begin{description}
        \item[\rm{(\mylabel{B.9}{$B.9$})}]
            There exists a constant $c > 0$ such that
            \begin{equation*}
                |\theta H(t, x, s \theta^{- 1}) - \zeta H(t, x, r \zeta^{- 1})|
                \leq c (1 + \max_{\tau \in [- h, t]} \|x(\tau)\|) \sqrt{\|s - r\|^2 + (\theta - \zeta)^2}
            \end{equation*}
            for any $t \in [0, T]$, $x(\cdot) \in C([- h, T], \mathbb{R}^n)$, $s$, $r \in \mathbb{R}^n$, and $\theta > 0$, $\zeta > 0$.

        \item[\rm{(\mylabel{B.10}{$B.10$})}]
            For every compact set $D \subset C([- h, T], \mathbb{R}^n)$, there exists a number $\lambda > 0$ such that
            \begin{equation*}
                \theta |H(t, x(\cdot), s \theta^{- 1}) - H(t, y(\cdot), s \theta^{- 1})|
                \leq \lambda \max_{\tau \in [- h, t]} \|x(\tau) - y(\tau)\|
            \end{equation*}
            for any $t \in [0, T]$, $x(\cdot)$, $y(\cdot) \in D$, and $s \in \mathbb{R}^n$, $\theta > 0$ such that $\|s\|^2 + \theta^2 = 1$.

        \item[\rm{(\mylabel{B.11}{$B.11$})}]
            For any $t \in [0, T]$, $x(\cdot) \in C([- h, T], \mathbb{R}^n)$, and $s \in \mathbb{R}^n$, the following limit exists:
            \begin{equation*}
                \lim_{\theta \downarrow 0} \theta H(t, x(\cdot), s \theta^{- 1})
                = H^0(t, x(\cdot), s),
            \end{equation*}
            and the functional $[0, T] \times C([- h, T], \mathbb{R}^n) \ni (t, x(\cdot)) \mapsto H^0(t, x(\cdot), s) \in \mathbb{R}$ is continuous for every fixed $s \in \mathbb{R}^n$.
    \end{description}

    Note that conditions \eqref{B.9} and \eqref{B.10} imply assumptions \eqref{B.2} and \eqref{B.3}, respectively.
    Note also that, for Hamiltonian \eqref{Hamiltonian_OCP}, assumptions \eqref{B.9} to \eqref{B.11} are fulfilled if the functions $f$ and $g$ satisfy conditions \eqref{A.1}, \eqref{A.3}, and the following condition instead of \eqref{A.2}:
    there exists a constant $c > 0$ such that
    \begin{equation*}
        \|f(\tau, y(\cdot), u)\| + |g(\tau, y(\cdot), u)|
        \leq c (1 + \max_{\xi \in [- h, \tau]} \|y(\xi)\|)
    \end{equation*}
    for every $\tau \in [0, T]$, $y(\cdot) \in C([- h, T], \mathbb{R}^n)$, and $u \in U$.

    Given $(\tau, y(\cdot)) \in [0, T] \times C([- h, T], \mathbb{R}^n)$ and $\overline{s} = (s, \theta) \in \mathbb{R}^n \times \mathbb{R}$, denote
    \begin{equation*}
        \begin{array}{c}
            \overline{E}^\ast(\tau, y(\cdot), \overline{s})
            = \big\{ (f, g) \in \overline{F}(\tau, y(\cdot)):
            \, \langle s, f \rangle + \theta g \geq \overline{H}(\tau, y(\cdot), \overline{s}) \big\}, \\[0.5em]
            \overline{E}_\ast(\tau, y(\cdot), \overline{s})
            = \big\{ (f, g) \in \overline{F}(\tau, y(\cdot)):
            \, \langle s, f \rangle + \theta g \leq \overline{H}(\tau, y(\cdot), \overline{s}) \big\}.
        \end{array}
    \end{equation*}
    Here,
    \begin{equation*}
        \begin{array}{c}
            \displaystyle
            \overline{F}(\tau, y(\cdot))
            = \big\{ (f, g) \in \mathbb{R}^n \times \mathbb{R}:
            \, \sqrt{\|f\|^2 + g^2} \leq \sqrt{2} c (1 + \max_{\xi \in [- h, \tau]} \|y(\xi)\|) \big\}, \\[1em]
            \overline{H}(t, x(\cdot), \overline{s})
            = \begin{cases}
                |\theta| H(t, x(\cdot), s |\theta|^{- 1}) & \mbox{if } \theta \neq 0, \\
                H^0(t, x(\cdot), s) & \mbox{if } \theta = 0,
            \end{cases}
        \end{array}
    \end{equation*}
    where the functional $H^0$ and the constant $c$ are taken from assumptions \eqref{B.9} and \eqref{B.10}.
    In this case, we have $(\mathbb{R}^n \times \mathbb{R}, \overline{E}^\ast) \in \mathcal{E}^\ast(H)$ and $(\mathbb{R}^n \times \mathbb{R}, \overline{E}_\ast) \in \mathcal{E}_\ast(H)$.
    Observe also that the Hamiltonian $\overline{H}$ is positively homogenous with respect to $\overline{s}$.

    By \cite[Theorem~1]{Krasovskii_Lukoyanov_2000_IMM_Eng}, there exists a unique functional $\varphi: [0, T] \times C([- h, T], \mathbb{R}^n) \to \mathbb{R}$ that is non-anticipative and continuous, meets \eqref{boundary_condition} and satisfies condition \eqref{MC} for $(\mathbb{R}^n \times \mathbb{R}, \overline{E}^\ast)$ and $(\mathbb{R}^n \times \mathbb{R}, \overline{E}_\ast)$.
    Then, according to Proposition~\ref{proposition_characteristic_complexes} and Theorem~\ref{theorem_existence_uniqueness}, this functional $\varphi$ is the minimax solution of problem \eqref{HJ}, \eqref{boundary_condition}.
    So, this paper also generalizes the results obtained for problem \eqref{HJ}, \eqref{boundary_condition} under assumptions \eqref{B.1} and \eqref{B.9} to \eqref{B.11}.

\section{Viscosity Solutions}
\label{section_viscosity}

    Following, e.g., \cite[Definition~2]{Crandall_Evans_Lions_1984}, a viscosity solution of problem \eqref{HJ}, \eqref{boundary_condition} can be naturally defined as a functional $\varphi: [0, T] \times C([- h, T], \mathbb{R}^n) \to \mathbb{R}$ that is non-anticipative and continuous, satisfies boundary condition \eqref{boundary_condition} and the following properties:
    \begin{description}
        \item[\rm{(\mylabel{V^ast}{$V^\ast$})}]
            For every $ci$-smooth functional $\psi: [0, T] \times C([- h, T], \mathbb{R}^n) \to \mathbb{R}$, if the difference $\varphi - \psi$ attains a local minimum at a point $(t, x(\cdot)) \in [0, T) \times C([- h, T], \mathbb{R}^n)$, then
    		\begin{equation*}
                \partial_t \psi(t, x(\cdot)) + H \big( t, x(\cdot), \nabla \psi(t, x(\cdot)) \big)
                \leq 0.
    		\end{equation*}
		\item[\rm{(\mylabel{V_ast}{$V_\ast$})}]
            For every $ci$-smooth functional $\psi: [0, T] \times C([- h, T], \mathbb{R}^n) \to \mathbb{R}$, if the difference $\varphi - \psi$ attains a local maximum at a point $(t, x(\cdot)) \in [0, T) \times C([- h, T], \mathbb{R}^n)$, then
    		\begin{equation*}
                \partial_t \psi(t, x(\cdot)) + H \big( t, x(\cdot), \nabla \psi(t, x(\cdot)) \big)
                \geq 0.
    		\end{equation*}
    \end{description}
    The statement below can be proved by the scheme from \cite[Theorem~1]{Lukoyanov_2007_IMM_Eng}.
    \begin{proposition} \label{proposition_minimax_viscosity}
        A minimax solution of problem \eqref{HJ}, \eqref{boundary_condition} is a viscosity solution of this problem.
    \end{proposition}

    Theorem~\ref{theorem_existence_uniqueness} and Proposition~\ref{proposition_minimax_viscosity} imply existence of a viscosity solution of problem \eqref{HJ}, \eqref{boundary_condition} under assumptions \eqref{B.1} to \eqref{B.3}.
    To obtain uniqueness results, a notion of a viscosity solution of path-dependent Hamil\-ton--Jacobi equations is usually modified by introducing auxiliary parameterizations (see, e.g., \cite{Soner_1988,Lukoyanov_2007_IMM_Eng,Lukoyanov_2010_IMM_Eng_1,Pham_Zhang_2014,Kaise_2015,Ekren_Touzi_Zhang_2016_1,Kaise_Kato_Takahashi_2018,Zhou_2019,Zhou_2020_1}), which allow to use compactness arguments.
    Uniqueness of viscosity solutions without such parameterizations is investigated in, e.g., \cite{Plaksin_2019_JOTA,Cosso_Russo_2019,Plaksin_2020_SIAM,Zhou_2020_2}.
    In \cite{Plaksin_2019_JOTA,Plaksin_2020_SIAM}, the problem is considered in the set of piecewise continuous functions and under more restrictive assumptions on the Hamiltonian and the boundary functional.
    In \cite{Cosso_Russo_2019,Zhou_2020_2}, suitable smooth variational principles are applied.
    Nevertheless, to the best of our knowledge, the question about uniqueness of viscosity solutions under assumptions \eqref{B.1} to \eqref{B.3} is open.

\end{document}